\documentclass [11pt] {article}

\usepackage {amssymb}

\title {Homogeneous Fedosov Star Products\\
        on Cotangent Bundles II: \\
        GNS Representations, the WKB Expansion, and Applications}

\author {{\bf
          Martin 
          Bordemann\thanks{Martin.Bordemann@physik.uni-freiburg.de}~,
          \addtocounter{footnote}{2}
          Nikolai 
          Neumaier\thanks{Nikolai.Neumaier@physik.uni-freiburg.de}~,
          Stefan 
          Waldmann\thanks{Stefan.Waldmann@physik.uni-freiburg.de}
         } \\[3mm]
         Fakult\"at f\"ur Physik\\Universit\"at Freiburg \\
         Hermann-Herder-Str. 3 \\
         79104 Freiburg i.~Br., F.~R.~G \\[3mm]
        }

\date{FR-THEP-97/23 \\[1mm]
      November 1997}

\newcommand {\LS} [1] {{#1} (\!(\lambda)\!)}
\newcommand {\NP} [1] {{#1} \langle\!\langle \lambda^* \rangle\!\rangle}
\newcommand {\CNP} [1] {{#1} \langle\!\langle \lambda \rangle\!\rangle}

\newcommand {\WV}   {\mathcal W \! \otimes \! \mbox{$\bigvee$}}
\newcommand {\WVL}  {\mbox{$\mathcal W \! \otimes \! \bigvee 
                           \! \otimes \Lambda$}}

\newcommand {\degs} {{\rm deg}_s}
\newcommand {\dega} {{\rm deg}_a}
\newcommand {\degh} {{\rm deg}_\lambda}
\newcommand {\Deg} {{\rm Deg}}

\newcommand {\TinyW} {{\mbox{\rm \tiny W}}}

\newcommand {\starw} {*_{\mbox{\rm \tiny W}}}

\newcommand {\wrep} {\varrho_{\TinyW}}

\newcommand {\TinyS} {{\mbox{\rm \tiny S}}}
\newcommand {\std} {\circ_\TinyS}
\newcommand {\ads} {{\rm ad}_\TinyS}
\newcommand {\rstd} {r_\TinyS}
\newcommand {\Dstd} {\mathcal D_\TinyS}
\newcommand {\taustd} {\tau_\TinyS}
\newcommand {\stars} {*_{\mbox{\rm \tiny S}}}
\newcommand {\srep} {\varrho_{\TinyS}}
\newcommand {\Srep} {\rho_\TinyS}
\newcommand {\vees} {{\star_\TinyS}}

\newcommand {\starf} {*_{\mbox{\rm \tiny F}}}

\newcommand {\im} {{\bf i}}
                                  
\newcommand {\BEQ} [1] {\begin {equation} \label {#1}}
\newcommand {\EEQ} {\end {equation}}

\newcommand {\Lie} {\mathcal L}
\newcommand {\cc} [1] {\overline {{#1}}}
\newcommand {\supp} {{\rm supp}}
\newcommand {\lsupp} {{\rm supp}_\lambda}
\newcommand {\id} {{\sf id}}
\newcommand {\tr} {{\sf tr}}
\newcommand {\ad} {{\rm ad}}
\newcommand {\SP} [2] {{\left\langle {#1}, {#2} \right\rangle}}
\newcommand {\Hom} {{\mathsf {Hom}}}
\newcommand {\End} {{\mathsf {End}}}
\newcommand {\field} [1] {{\mathsf {{#1}}}}
\newcommand {\Div} {{\mathsf {div}}}
\newcommand {\graph} {{\rm graph}}

\newcommand {\Hfree} {{H_{\mbox{\tiny free}}}}

\newenvironment {PROOF}{\small {\sc Proof:}}{{\hspace*{\fill} $\square$}}

\newtheorem {LEMMA} {Lemma} [section]
\newtheorem {PROPOSITION} [LEMMA] {Proposition}
\newtheorem {THEOREM} [LEMMA] {Theorem}
\newtheorem {COROLLARY} [LEMMA] {Corollary}
\newtheorem {DEFINITION}[LEMMA] {Definition}

\begin {document}

\maketitle

\begin {abstract}
   
This paper is part II of a series of papers on the deformation
quantization on the cotangent bundle of an arbitrary manifold $Q$.
For certain homogeneous star products of Weyl ordered type (which we 
have obtained from a Fedosov type procedure in part I, see q-alg/9707030) 
we construct differential operator representations via the formal GNS 
construction (see q-alg/9607019). The positive linear functional is 
integration over $Q$ with respect to some fixed density and is shown 
to yield a reasonable version of the Schr\"odinger representation
where a Weyl ordering prescription is incorporated. 
Furthermore we discuss simple examples like free particle Hamiltonians 
(defined by a Riemannian metric on $Q$) and the implementation of certain 
diffeomorphisms of $Q$ to unitary transformations in the GNS (pre-)Hilbert 
space and of time reversal maps (involutive anti-symplectic 
diffeomorphisms of $T^*Q$) to anti-unitary transformations. We show that 
the fixed-point set of any involutive time reversal map is either empty 
or a Lagrangean submanifold. Moreover, we compare our approach to
concepts using integral formulas of generalized Moyal-Weyl type. 
Furthermore we show that the usual WKB 
expansion with respect to a projectable Lagrangean submanifold can be 
formulated by a GNS construction. Finally we prove that any homogeneous 
star product on any cotangent bundle is strongly closed, i.~e. the integral 
over $T^*Q$
w.r.t.~the symplectic volume vanishes on star-commutators. An alternative 
Fedosov type deduction of the star product of standard ordered type 
using a deformation of the algebra of symmetric contravariant tensor fields 
is given.

\end {abstract}

\newpage


\section {Introduction}
\label {IntroSec}

Deformation quantization is a now well-established quantization 
concept introduced by Bayen, Flato, Fronsdal, Lichnerowicz, and 
Sternheimer \cite{BFFLS78}. For any symplectic manifold
the existence of the formal associative deformation (the star product)
of the pointwise multiplication of smooth functions which are identified 
with the classical observables had been shown by \cite{DL83} and 
\cite{Fed94}. A recent preprint by Kontsevich states that such star 
products even exist for any Poisson manifold \cite {Kon97}. 
The classification for star products on symplectic manifolds 
up to equivalence transformations by formal power series
in the second de Rham cohomology group is due to \cite{NT95a,NT95b} and
\cite{BCG96}.

The particular case of a cotangent bundle $T^*Q$ where $Q$ is the 
physical configuration space is of great importance for physicists and 
there is a large amount of literature concerning various ways of 
quantization of such physical phase spaces. 
Deformation quantization and star products on cotangent bundles are 
considered e.~g. in \cite {CFS92,DL83a,Pfl95}.
Moreover differential operator representations and their symbolic 
calculus including integration techniques are considered in e.~g.
\cite{Und78,Wid80,Emm93a, Emm93b}. For geometric quantization of 
cotangent bundles see e.~g. \cite{Woo80} and references therein.

In order to formulate states and Hilbert space representations of the 
deformed algebra of formal power series of smooth complex-valued functions 
on a symplectic manifold two of us (M.~B. and S.~W.) have recently 
transferred the concept of GNS representations in the theory of 
$C^*$-algebras to deformation quantization \cite {BW96b}. In that paper 
this concept was shown to be physically reasonable by applying it mainly 
to flat $\mathbb R^{2n}$: integration over configuration space at fixed 
momentum value $0$ turned out to be a formally positive linear functional 
whose GNS representation yielded the usual formal Schr\"odinger 
representation as formal differential  operators on $\mathbb R^n$ by means 
of the Weyl ordering prescription (see \cite {BNW97a} for details). In 
another paper these results had been extended to the following situation:
integration over a projectable submanifold $L$ of $\mathbb R^{2n}$ with 
respect to a geometrically defined volume preceded by a geometrically 
defined formal series of differential operators yielded a formally 
positive linear functional whose GNS representation contained the usual 
WKB expansion of an eigenfunction of some Hamiltonian operator to its 
eigenvalue $E$ provided $L$ is contained in the energy surface of the 
classical Hamiltonian to the value $E$ \cite {BW97b}.

This paper is part II of a series of papers devoted to the study of 
certain star products on the cotangent bundle of an arbitrary manifold $Q$ 
and applications of the GNS construction for this particular physically 
important situation. In part I (see \cite {BNW97a}) we had prepared the 
ground by constructing a star product of standard ordered type $\stars$ 
(in which the first function is differentiated in vertical 
`momentum' directions only) and a representation of this algebra where 
those functions on $T^*Q$ canonically corresponding to symmetric 
contravariant tensor fields $T$ on $Q$ (i.~e. which are polynomial 
in the momenta) are mapped to differential operators on $Q$ where 
$T$ is naturally paired 
with a multiple covariant derivative with respect to some torsion-free 
connection on $Q$. This had been done by a Fedosov type construction. 
Moreover we had constructed a star product $\starw$ of Weyl type on 
$T^*Q$ equivalent to $\stars$ and a ${ }^*$-representation (with respect to 
complex conjugation) of this algebra as generalized Weyl ordered 
differential operators.

The aim of this paper is to generalize the results of \cite {BW96b,BW97b} 
to arbitrary cotangent bundles, that is to obtain the Schr\"odinger-like 
differential operator representation as well as the WKB expansion as GNS 
representations. Moreover we discuss some simple applications and 
compare our results with the approaches using integral 
formulas like for instance Underhill's quantization \cite {Und78}
and \cite {Pfl97,Emm93a,Emm93b} and give an elementary proof that any 
homogeneous star product on any cotangent bundle is strongly closed in the 
sense of \cite {CFS92}.

The paper is organized as follows: 
Section \ref {PrelimSec} contains notation and results 
of the first part. In section \ref {FedStarsSec} we give 
another independent construction for the star product $\stars$ which 
provides another way to show that this star product is differential. In 
section \ref {GNSSec} we consider a linear functional which is 
integration over the configuration space $Q$ with respect to some chosen 
positive density and show that this turns out to be a formally positive 
linear functional. Hence it induces 
a GNS representation which is explicitly computed. In the following 
section \ref {ApplSec} we mention some easy but physically important
convergence properties and discuss simple physical examples. 
Moreover we consider particular
symmetries of $T^*Q$ and their implementation as automorphisms of 
the star product algebras and as unitary maps in the GNS Hilbert spaces.
We prove that involutive anti-symplectic diffeomorphisms have either no 
fixed points or the fixed-point set is a Lagrangean submanifold and
describe such maps for certain coadjoint orbits.
In section \ref {CompSec} we compare our formulas using integral formulas 
for the representations with other approaches already mentioned. 
Section \ref {WKBSec} is devoted to the WKB expansion where we 
generalize the previous results to the case of projectable Lagrangean 
submanifolds in a cotangent bundle. Finally we prove in section 
\ref {TraceSec} that any homogeneous star product on $T^*Q$ is strongly 
closed.

\section {Preliminary results and notation}
\label {PrelimSec}

In this section we shall remember some results from \cite {BW96b,BNW97a} 
and establish our notation:
Let $Q$ be a smooth $n$-dimensional manifold and $\pi: T^*Q \to Q$ its 
cotangent bundle with the canonical one-form $\theta_0$ and the
canonical symplectic form $\omega_0 = -d\theta_0$. 
Moreover we denote by $\xi$ the canonical Liouville vector field 
defined by $i_\xi \omega_0 = - \theta_0$. 
Let $i: Q \to T^*Q$ be the canonical embedding of $Q$ as zero section. 
Moreover we shall fix a torsion-free connection $\nabla$ on $Q$ and 
make use of local bundle Darboux coordinates
$q^1, \ldots, q^n, p_1, \ldots, p_n$ induced by local coordinates 
$q^1, \ldots, q^n$ on $Q$. A function $f \in C^\infty (T^*Q)$ is called 
polynomial in the momenta of degree $k$ if $\Lie_\xi f = kf$ and 
the set of polynomial functions in the momenta of degree $k$ is denoted by 
$C^\infty_{pp,k} (T^*Q)$ and we set 
$C^\infty_{pp} (T^*Q) := \bigoplus_{k=0}^\infty C^\infty_{pp,k} (T^*Q)$.
Then $C^\infty_{pp} (T^*Q)$ is canonically isomorphic as graded algebra 
to $\Gamma (\bigvee TQ)$ together with the symmetric product $\vee$ 
via the canonical isomorphism 
$\Gamma (\bigvee^k TQ) \ni T \mapsto \widehat T \in C^\infty_{pp,k} (T^*Q)$
given by 
$\widehat T (\alpha_q) := \frac{1}{k!} T(\alpha_q, \ldots, \alpha_q)$
where $\alpha_q \in T^*_q Q$. For some easy homogeneity properties of 
$C^\infty_{pp} (T^*Q)$ see e.~g. \cite [Lemma A.2] {BNW97a}.
We shall use Einstein's summation convention, i.~e. summation over 
repeated indices is understood.

Throughout this paper we denote the formal deformation parameter by 
$\lambda$ which corresponds directly to Planck's constant $\hbar$.
Using the connection $\nabla$ on $Q$ we had constructed a star product of 
standard ordered type $\stars$ in \cite {BNW97a} which is the 
natural generalization of the standard ordered product in flat 
$\mathbb R^{2n}$. 
A more physically star product of Weyl type having the complex conjugation 
as involutive antilinear anti-automorhism has been shown to exist and is 
constructed by the following equivalence transformation: 
we consider the operator 
\BEQ {DefDelta}
  \Delta:= \partial_{p_i} \partial_{q^i} 
           + p_r (\pi^* \Gamma^r_{ij}) \partial_{p_i} \partial_{p_j}
           + (\pi^*\Gamma^i_{ij}) \partial_{p_j}
           + (\pi^* \alpha_j) \partial_{p_j}
\EEQ
where locally $\alpha = \alpha_j dq^j$ is a one-form satisfying 
$\tr R = -d\alpha$ where $R$ is the curvature tensor of $\nabla$. This 
operator is globally defined and induces the equivalence transformation
(see \cite [Section 7] {BNW97a})
\BEQ {DefNeumaier}
        N:= e^{\frac{\lambda}{2\im}\Delta}
\EEQ
which yields the star product $\starw$ by 
\BEQ {DefStarw}
         f\starw g := N^{-1}((Nf)\stars(Ng))
\EEQ
where $f, g \in C^\infty (T^*Q)[[\lambda]]$. 
This equivalence is again the natural generalization of the flat case and 
hence we shall call $\starw$ the star product of Weyl ordered type.
Both star products have the important property of being homogeneous 
in the sense of \cite {DL83a}: the following operator
\BEQ {HomoDef}
    \mathcal H := \lambda\frac{\partial}{\partial \lambda} + \Lie_\xi
\EEQ 
is a derivation of $\stars$ and $\starw$. Furthermore a representation on 
$C^\infty (Q)[[\lambda]]$ for the standard ordered product 
$\stars$ was shown to be given by the following explicit formula
\BEQ {srepDef}
    \srep (f) \psi := i^* (f \stars \pi^*\psi) =
    \sum_{r=0}^\infty \frac{1}{r!}
    \left(\frac{\lambda}{\im}\right)^r
    i^*\left( \frac{\partial^r f}{\partial p_{i_1} \cdots
    \partial p_{i_r}} \right)
    i_s (\partial_{q^{i_1}}) \cdots i_s (\partial_{q^{i_r}})
    \frac{1}{r!} D^r \psi
\EEQ
where $f \in C^\infty (T^*Q)[[\lambda]]$ and 
$\psi \in C^\infty (Q)[[\lambda]]$ and where 
$D := dq^k \vee {\nabla}_{\partial q^k}$ justifying the notion 
of `standard ordered type' by analogy to the flat case. 
Then a representation $\wrep$ for $\starw$ on the 
same representation space is given by
\BEQ {wrepDef}
     \wrep (f)\psi := \srep (Nf) \psi = i^*((Nf) \stars \pi^* \psi)     
\EEQ     
yielding the explicit formula
\BEQ {wrepExpl}
    \wrep (f) \psi = \sum_{r=0}^\infty \frac{1}{r!}
    \left(\frac{\lambda}{\im}\right)^r
    i^*\left( \frac{\partial^r Nf}{\partial p_{i_1} \cdots
    \partial p_{i_r}} \right)
    i_s (\partial_{q^{i_1}}) \cdots i_s (\partial_{q^{i_r}})
    \frac{1}{r!} D^r \psi .
\EEQ

In \cite {BW96b} a generalization of the well-known GNS construction for 
$C^*$-algebras to deformation quantization was proposed where the crucial 
ingredient is the notion of positivity in the ring of formal power series 
$\mathbb R[[\lambda]]$ (see e.~g. appendix \ref {CNPApp}). We have argued 
that it is advantageous to consider the field of formal Laurent series 
instead of $\mathbb R[[\lambda]]$ which is in a canonical way an ordered 
field. Moreover it turned out that certain field extensions (the formal 
NP and CNP series) are even more suitable for the definition of the 
GNS representation and `formal' Hilbert spaces. In appendix \ref {CNPApp} 
we provide several lemmata and definitions how the results obtained for 
the more convenient formal power series can be extended to these more 
general formal series justifying thereby our sloppy usage of the notion 
of `Hilbert' spaces etc. In fact all relevant structures are completely 
determined in the setting of formal power series, see lemma 
\ref {ExtendLem} and \ref {ExtendPosFunctLem}.

\section {Another Fedosov-like construction of the standard ordered 
          star product $\stars$}
\label {FedStarsSec}

In this section we shall give a different method to obtain the 
standard ordered star product $\stars$ as in \cite {BNW97a} 
avoiding the additional technical ingredient of a lifted connection. 
To do so we first notice that it is enough to consider
$C^\infty_{pp}(T^*Q)[[\lambda]]$ instead of
$C^\infty(T^*Q)[[\lambda]]$ since the star product is given by
bidifferential operators, which are uniquely determined by their
values on $C^\infty_{pp}(T^*Q)$. Using the canonical algebra 
isomorphism $\widehat{\quad}$ between 
$\Gamma(\bigvee TQ)$ and $C^\infty_{pp}(T^*Q)$ 
the task of finding a deformation of the pointwise multiplication in
$C^\infty_{pp}(T^*Q)$ is equivalent to the deformation of the
$\vee$-product in the symmetric algebra $\Gamma(\bigvee TQ)$.
Therefore we start a Fedosov procedure similar to the 
construction in \cite {Bor96} to deform the $\vee$-product
using a connection on $Q$ and show that this deformation is
compatible to the standard ordered representation of functions that
are polynomial in the momentum variables, giving an alternative
method to construct the star product $\stars$.

\subsection*{A deformation $\vees$ of the $\vee$-product}

We consider the slightly modified Fedosov algebra
\[
    \WVL := \left( {\mathsf X}_{s=0}^\infty
           \mathbb C \left(\Gamma\left(\mbox{$\bigvee$}^s T^*Q \otimes   
           \mbox{$\bigvee$} TQ \otimes \mbox{$\bigwedge$}
           T^*Q \right)\right)\right)[[\lambda]].
\]
To avoid clumsy notation we drop the explicit mention of the
complexification that shall be taken for granted. To define the
important mappings we make use of factorized sections of the shape
$F_i = \lambda^{q_i} f_i \otimes S_i \otimes \alpha_i$ with 
$f_i \in \Gamma(\bigvee^{s_i} T^*Q)$, $S_i \in \Gamma(\bigvee^{d_i} TQ)$,
$\alpha_i \in \Gamma(\bigwedge^{a_i} T^*Q)$. In addition we
consider the obvious degree maps $\degs$, $\degs^*$,
$\dega$, and $\degh$ with
\[
    \degs F_i= s_i F_i, 
    \quad \degs^* F_i = d_i F_i,
    \quad \dega F_i = a_i F_i, 
    \quad \degh F_i = q_i F_i .
\]
The additional symmetric degree with respect to $\degs^*$
is referred to as {\em dual} symmetric degree. The undeformed
multiplication in $\WVL$ of two sections $F_1, F_2$ is defined by
\[
    \mu \circ (F_1 \otimes F_2) = F_1 F_2 := 
    \lambda^{q_1 + q_2} (f_1 \vee f_2) \otimes (S_1 \vee S_2) 
    \otimes (\alpha_1 \wedge \alpha_2) .
\]
By $\sigma$ we denote the linear map 
$\sigma : \WVL \to \Gamma (\bigvee TQ)[[\lambda]]$ that projects onto 
the part of symmetric and antisymmetric degree zero. Given a torsion-free
connection $\nabla$ on $Q$ we define the connection $\nabla$ in
$\WVL$ using a chart $q^1, \ldots, q^n$ of $Q$ by
\[
    \nabla := ( 1 \otimes 1 \otimes dq^l ) \nabla_{\partial_{q^l}},
\]
which turns out to be a globally defined homogeneous super derivation 
of the pointwise product of degree $(0,0,1)$. Now we define the fibrewise
associative deformed multiplication $\std$ for $F, G \in \WVL$ by
\[
    F \std G := \mu \circ 
                e^{\frac{\lambda}{\im} 
                i_s^*(dq^l) \otimes i_s(\partial_{q^l})} F \otimes G,
\]
where $i_s^*(dq^l)$ denotes the symmetric insertion of the one-form
$dq^l$ and $i_s(\partial_l)$ denotes the symmetric insertion of the
vector field $\partial_{q^l}$. Following Fedosov we define the mappings
\[
    \delta := (1 \otimes 1 \otimes dq^l) i_s(\partial_{q^l})
    \quad \mbox { and } \quad 
    \delta^* := (dq^l \otimes 1 \otimes 1) i_a(\partial_{q^l}),
\]
which are super derivations of the undeformed product of degree
$(-1,0,1)$ and $(1,0,-1)$ and for which the following identities
hold:
\[
    \delta^2 = {\delta^*}^2 = 0, 
    \quad 
    \delta \delta^* + \delta^* \delta = \degs + \dega, 
    \quad 
    \nabla \delta + \delta \nabla = 0.
\]
In addition one easily verifies, that $\delta$ resp. $\nabla$ are
super derivations with respect to $\std$ of degree
$(-1,0,1)$ resp. $(0,0,1)$. In analogy to the known Fedosov
construction we define
\[
    \delta^{-1} F_1 := 
    \left\{ \begin{array}{cl}
            \frac{1}{s_1+a_1} \delta^* F_1 
            & \textrm{ for } s_1+a_1 \neq 0 \\
            0 
            & \textrm{ for } s_1+a_1 = 0
            \end{array} 
    \right.
\]
with which a `Hodge-decomposition' for any $F \in \WVL$ is valid, i.~e.
\[
    F = \delta \delta^{-1}F + \delta^{-1} \delta F + \sigma(F).
\]
Finally we define $\dega$-graded super commutators with
respect to $\std$ by 
$[F_1, F_2] := \ad_\TinyS(F_1)F_2 := 
F_1 \std F_2 - (-1)^{a_1a_2} F_2 \std F_1$ 
and the total degree map $\Deg := 2 \degh + \degs + \degs^*$ 
being a $\std$ super derivation of degree $(0,0,0)$. 
The operator corresponding to the homogeneity derivation
$\mathcal H$ as in (\ref{HomoDef}) is given by 
$\mathsf H := \degs^* +\degh $ and is a (only $\mathbb C$-linear!) 
derivation of $\std$. 
After these preparations we have the first little lemma being of
importance for the construction of the Fedosov derivation:
\begin{LEMMA}
For $R_\TinyS := -\frac{1}{2} R^l_{kij} dq^k \otimes \partial_{q^l}
\otimes dq^i \wedge dq^j$ denoting by $R^l_{kij}$ the components
of the curvature tensor of the connection $\nabla$ we have
\[
    \nabla^2 = \frac{\im}{\lambda} \ads (R_\TinyS), 
    \qquad 
    \delta R_\TinyS =0,
    \qquad 
    \nabla R_\TinyS =0.
\]
\end{LEMMA}
\begin{PROOF}
The first assertion is a straight forward computation and the
others are just reformulations of the Bianchi-identities.
\end{PROOF}

For the Fedosov derivation $\Dstd$ satisfying $\Dstd^2=0$ we make the 
usual Ansatz
\[
    \Dstd = -\delta + \nabla + \frac{\im}{\lambda} \ads (\rstd)
\]
with $\rstd \in \WVL^1$ (i.~e. $\dega \rstd = \rstd$) and completely 
analogously to the usual Fedosov procedure 
(see e.~g. \cite {Fed94,BNW97a}) we obtain
a uniquely determined such $\rstd$ obeying the conditions 
$\delta^{-1}\rstd=0$ and 
$\delta \rstd = 
\nabla \rstd + R_\TinyS + \frac{\im}{\lambda} \rstd \std \rstd$. 
Moreover $\rstd = \sum_{k=3}^\infty \rstd^{(k)}$ with 
$\Deg \rstd^{(k)} = k \rstd$ is given by the following 
recursion formula
\[
    \rstd^{(3)} = \delta^{-1} R_\TinyS \qquad 
    \rstd^{(k+3)} = \delta^{-1} \left( 
                    \nabla \rstd^{(k+2)} + 
                    \frac{\im}{\lambda} \sum_{l=1}^{k-1}
                    \rstd^{(l+2)} \std \rstd^{(k-l+2)}
                    \right).
\]
Using this recursion formula one can show by induction on the total degree 
that $\mathsf H \rstd = \rstd$ and hence $\mathsf H \Dstd = \Dstd \mathsf H$
implying together with the particular shape of $\std$ the following easier 
recursion formula
\[
    \rstd^{(3)} = \delta^{-1} R_\TinyS 
    \qquad
    \rstd^{(k+3)} = \delta^{-1}\left( 
                    \nabla \rstd^{(k+2)} 
                    - \frac{1}{2} \sum_{l=1}^{k-1} 
                    \left\{\rstd^{(l+2)} , \rstd^{(k-l+2)}
                    \right\}_{\mathrm{fib}}
                    \right)
\]
denoting by $\{{}\cdot{},{}\cdot{}\}_{\mathrm{fib}}$ the fibrewise
Poisson bracket given by 
$\{F,G\}_{\mathrm{fib}} = i_s(\partial_{q^l}) F i_s^*(dq^l) G -
i_s^*(dq^l) F i_s(\partial_{q^l})G$. 
Obviously this implies that $\rstd$ does not depend on $\lambda$ at all.
Moreover we have the following proposition:
\begin{PROPOSITION}
$\mathcal{W}_{\Dstd} := ker( \Dstd) \cap \WV$ is a subalgebra 
of $(\WVL, \std)$ and the map $\sigma$ restricted to $\mathcal{W}_{\Dstd}$ 
is a $\mathbb C[[\lambda]]$-linear bijection onto 
$\Gamma(\bigvee TQ)[[\lambda]]$. Denoting by 
$\taustd: \Gamma(\bigvee TQ)[[\lambda]] \to 
\mathcal{W}_{\Dstd} \subset \WV$ 
the inverse of the restriction of $\sigma$ to
$\mathcal{W}_{\Dstd}$ we have the following recursion
scheme to calculate $\taustd(T)$ for $T = \sum_{t=0}^m T^{(t)}$ 
$(m \in \mathbb N)$ and $T^{(t)} \in \Gamma(\bigvee^t TQ)$:
The Fedosov-Taylor series 
$\taustd(T) = \sum_{l=0}^\infty \taustd (T)^{(l)}$ with 
$\Deg \taustd (T)^{(l)}= l \taustd(T)^{(l)}$ is given by
\begin{eqnarray} \label{taurec}
    \taustd(T)^{(0)}   
    & = & T^{(0)} \\ 
    \taustd(T)^{(k+1)} 
    & = & \delta^{-1} \left( 
          \nabla \taustd(T)^{(k)} + \frac{\im}{\lambda}
          \sum_{t=1}^{k-1} \ads \left(\rstd^{(t+2)} \right) 
          \taustd(T)^{(k-t)} \right) 
          + T^{(k+1)} \mbox{ for } 0 \leq k \leq m-1 \\
    \taustd(T)^{(k+1)} 
    & = & \delta^{-1} \left( 
          \nabla \taustd(T)^{(k)} + \frac{\im}{\lambda}
          \sum_{t=1}^{k-1} \ads \left(\rstd^{(t+2)} \right) 
          \taustd(T)^{(k-t)} \right) 
          \mbox{ for } k \geq m .
\end{eqnarray}
In addition the $\mathbb C [[\lambda]]$-linear mapping $\taustd$ satisfies
\BEQ {TauHomogen}
    \mathsf H \taustd(T) = \taustd(\mathsf H T) 
    \quad \forall T \in \Gamma(\mbox{$\bigvee$} TQ)[[\lambda]].
\EEQ
\end{PROPOSITION}
\begin{PROOF}
The recursion formula is proven analogously to 
\cite [Theorem 2.2] {BNW97a} and the homogeneity of $\taustd$ follows 
directly from $\mathsf H \Dstd = \Dstd \mathsf H$.
\end{PROOF}

Remark: Observing the homogeneity of $\taustd$ with respect to
$\mathsf H$ the recursion formula for 
$\taustd(T)=\sum_{k=0}^\infty \taustd(T)_{(k)}$ with 
$(\degs + \degh) \taustd(T)_{(k)} = k\taustd(T)_{(k)}$ for 
$T \in \Gamma(\bigvee TQ)$ can be rewritten in a more convenient 
manner, i.~e.
\[
    \taustd(T)_{(0)} = T,
    \qquad
    \taustd(T)_{(k+1)} = \delta^{-1} \left(
                         \nabla \taustd(T)_{(k)} + 
                         \frac{\im}{\lambda} \sum_{t=1}^k 
                         \ads \left(\rstd^{(t+2)} \right) 
                         \taustd(T)_{(k-t)} \right).
\]
Now the associative product $\vees$ for $\Gamma(\bigvee TQ)[[\lambda]]$ 
is defined by pull-back of $\std$ via $\taustd$, i.~e.
\[
    S \vees T := \sigma( \taustd (S) \std \taustd(T))
\]
for $S, T \in \Gamma(\bigvee TQ)[[\lambda]]$, that induces an
associative product $\stars$ for $C^\infty_{pp}(T^*Q)[[\lambda]]$
by pull-back via the natural isomorphism $\widehat{\quad}$ i.~e.
\BEQ {Alterstars}
    \widehat{S} \stars \widehat{T} := \widehat{(S \vees T)}
\EEQ
for $\widehat{S},\widehat{T} \in C^\infty_{pp}(T^*Q)[[\lambda]]$. 
We shall now prove that this definition indeed coincides with the standard 
ordered product as in \cite {BNW97a}. 
Using (\ref {TauHomogen}) and the fact that $\mathsf H$ is a 
$\vees$-derivation one can show that $\mathcal H$ is a derivation 
of $\stars$ and thus $\stars$ is homogeneous.

\subsection* {A representation of $(\Gamma(\bigvee TQ)[[\lambda]],
              \vees)$}

This subsection is dedicated to construct a representation of
$(\Gamma(\bigvee TQ)[[\lambda]], \vees)$, that is compatible with the 
standard representation we constructed in \cite [Section 6] {BNW97a}. 
As a starting point we find a fibrewise representation of 
$(\WV , \std)$. First we define the representation space
\[
    \mathfrak H := \WV \cap \ker (\degs^*)
\]
and the projection 
$P : \WV \to \WV \cap \ker (\degs^*)$ that projects onto 
the subspace of dual symmetric degree zero. For $T \in \WV$ and 
$\Psi \in \mathfrak H$ we define the fibrewise standard representation
$\widetilde{\Srep}: \WV \to \End (\mathfrak H)$ by
\BEQ {fiberrep}
    \widetilde{\Srep}(T) \Psi := P(T \std \Psi).
\EEQ
\begin{LEMMA}
The map $\widetilde{\Srep}$ is a $\std$-representation of $\WV$
on $\mathfrak H$, i.~e. for $S, T \in \WV$ we have
\[
    \widetilde{\Srep}(S \std T) = 
    \widetilde{\Srep}(S) \widetilde{\Srep}(T).
\]
\end{LEMMA}
\begin{PROOF}
The $\mathbb C [[\lambda]]$-lineartity of $\widetilde{\Srep}$ is
obvious and the representation property is proved by straight
forward computation using the associativity of $\std$ and the
validity of the equation $P(F \std G) = P(F \std P G)$ for 
$F, G \in \WV$, which follows from the particular shape of $\std$.
\end{PROOF}

What we have in mind is to define a $\vees$-representation on the
representation space $C^\infty(Q)[[\lambda]]$ that can canonically
be identified with $\WV \cap \ker(\degs^*)\cap \ker(\degs)$, 
thus for $\psi \in C^\infty(Q)[[\lambda]]$ and 
$T \in \Gamma(\bigvee TQ) [[\lambda]]$ the following expression is
well-defined
\BEQ {rep}
    \Srep(T)\psi := p(T \vees \psi),
\EEQ
denoting by $p$ the projection from $\Gamma(\bigvee TQ) [[\lambda]]$ 
to the part of dual symmetric degree zero, 
i.~e. $C^\infty (Q) [[\lambda]]$. The representation property 
of $\Srep$ can now be easily proven:
\begin{PROPOSITION}
Let $\Dstd$ be the Fedosov derivation constructed
above and let $\taustd$ be the corresponding Fedosov-Taylor series
and $\vees$ the deformation of the $\vee$-product. Then
\BEQ {PcommD}
    \Dstd P = P \Dstd
\EEQ
and
\BEQ {stdrep}
    \Srep(T)\psi = p(T \vees \psi) = 
    \sigma( \widetilde{\Srep} (\taustd(T))\taustd(\psi))
\EEQ
defines a $\vees$-representation of $\Gamma (\bigvee TQ)[[\lambda]]$ 
where $\psi \in C^\infty(Q)[[\lambda]]$ and 
$T \in \Gamma(\bigvee TQ)[[\lambda]]$. 
\end{PROPOSITION}
\begin{PROOF}
As an abuse of notation we also denoted by $P$ the mapping
$P\otimes 1 : \WVL^1 \to \WVL^1 \cap \ker(\deg_s^*)$.
Then one immediately checks that $\delta$ and $\nabla$ commute with $P$.
Moreover $\ad_\TinyS(\rstd)$ commutes with $P$ due to the
particular shape of $\std$ and $\degs^* \rstd =
\rstd$. Using this and the obvious equation $p (\sigma(F))=
\sigma (P(F))$ for $F \in \mathcal W$ it is straight forward
proving the representation property of $\Srep$ observing equation
(\ref{PcommD}) and $P(F \std G) = P(F \std P G)$
for $F, G \in \WV$.
\end{PROOF}

To conclude this section we shall now show the compatibility of the
representation $\srep$ to the standard representation constructed
in \cite{BNW97a}, implying that $\stars$ defined on
$C^\infty_{pp}(T^*Q)[[\lambda]]$ by pull-back of $\vees$ as above
coincides with the standard ordered star product we considered in
\cite{BNW97a}. Moreover the star product $\stars$ for 
$f, g \in C^\infty(T^*Q)[[\lambda]]$ is completely determined by this
construction since bidifferential operators are completely
determined by their values on functions which are polynomial in the  
momenta.
\begin{PROPOSITION}
For $\psi \in C^\infty(Q)[[\lambda]]$ the Fedosov-Taylor series is
given by
\BEQ {Taylorfunct}
    \taustd ( \psi ) = e^D \psi \quad 
    \textrm{ with} \quad D := dq^l \vee \nabla_{\partial_{q^l}} .
\EEQ
Therefore the representation $\Srep$ can be calculated explicitly
for $T \in \Gamma(\bigvee TQ)[[\lambda]]$ and is given by
\BEQ {Standrep}
    \Srep(T) \psi = \sum_{r=0}^\infty \frac{1}{r!} 
                    \left(\frac{\lambda}{\im} \right)^r 
                    p\left( i_s^*(dq^{i_1}) \ldots i_s^*(dq^{i_r})T \right) 
                    i_s(\partial_{q^{i_1}}) \ldots i_s(\partial_{q^{i_r}})
                    \frac{1}{r!} D^r \psi 
\EEQ
and thus
\BEQ {Srepsrep}
    \Srep (T) \psi = \srep (\widehat T) \psi .
\EEQ    
Hence the obvious injectivity of $\Srep$ implies that the product defined 
by (\ref {Alterstars}) indeed coincides with $\stars$.
\end{PROPOSITION}
\begin{PROOF}
The formula for $\taustd (\psi)$ can be proven by induction on the
total degree using the recursion formula given in (\ref{taurec}).
The equation (\ref{Standrep}) then follows by direct calculation
from the definition of the representation $\Srep$ in equation
(\ref{rep}). Clearly (\ref {Srepsrep}) follows by direct 
comparison with (\ref {srepDef}).
\end{PROOF}

At last in this section we should mention that our considerations yield an 
alternative proof of the fact that the standard representation for functions
polynomial in the momenta gives rise to a standard ordered star product 
$\stars$ of Vey type since we have the following proposition:
\begin{PROPOSITION}
For $S, T \in \Gamma(\bigvee TQ)$ there are bidifferential operators 
$\mathsf{N}_r$ $(r \in \mathbb N)$ on $\Gamma(\bigvee TQ)$ that are of 
order $r$ in both arguments such that $S\vees T$ can be written as 
\[
    S \vees T = \sum_{r=0}^\infty
                \left(\frac{\lambda}{\im}\right)^r 
                \mathsf N_r(S,T).
\]
\end{PROPOSITION}
\begin{PROOF}
The proof is an easy consequence of the definition of $\vees$ and the fact 
that the mapping $\Gamma(\bigvee TQ) \ni T \mapsto  \taustd(T)_{(r),l} \in 
\Gamma(\bigvee^l T^*Q \otimes \bigvee TQ)$ is a differential operator of 
order $r$ for all $0\leq l \leq r$ (where we have written 
$\taustd(T)_{(r)} = \sum_{l=0}^r  \left(\frac{\lambda}{\im}\right)^{r-l} 
\taustd(T)_{(r),l}$), what can be proven by straight forward induction 
on $r$ using the recursion formula for $\taustd(T)_{(r)}$.
\end{PROOF}\\
Remark: At this instance we should mention that we use the notion of 
differential operators on the commutative, associative algebra 
$(\Gamma(\bigvee TQ),\vee)$ as a purely algebraic concept. 
For the reader unfamiliar with this point of view, we just give a short 
definition: By $l_T$ we denote the left multiplication by 
$T\in \Gamma(\bigvee TQ)$ that is defined to be a differential operator 
of order $0$. Recursively the differential operators of order 
$k \in \mathbb N$ are defined by the set of all 
endomorphisms $\mathcal D$ of $\Gamma(\bigvee TQ)$ satisfying that 
$[\mathcal D, l_T]$ is a differential operator of order $k-1$. 
Similarly one defines multi-differential operators between modules over a 
commutative and associative algebra. In view of 
this definition the insertions $i_s^*(\beta)$ for a one-form $\beta$ have 
to be considered as differential operators of order $1$ on 
$\Gamma(\bigvee TQ)$ since for all $T \in \Gamma(\bigvee TQ)$ we have 
$[i_s^*(\beta), l_T]=l_{i_s^*(\beta)T}$.

Moreover the conjugation with the natural algebra-isomorphism 
$\widehat{\quad}$ leaving invariant the order of differential 
operators yields bidifferential operators $N_r$ of order $r$ in 
every argument on $C_{pp}^\infty(T^*Q)$ given by 
$N_r(\widehat S, \widehat T) = \widehat{\mathsf N_r(S,T)}$ 
implying that the star product $\stars$ is of Vey type.

\section {A GNS construction for a Schr\"odinger representation} 
\label {GNSSec}

For a given complex number $\gamma$ let
$|\!\bigwedge\!|^\gamma T^*Q$ be the bundle of $\gamma$-densities on $Q$
(see e.~g.~\cite[p.119--121]{BW95}): this bundle can be obtained by
taking the bundle of linear frames $L(Q)$ over $Q$ and by associating
the typical fibre $\mathbb C$ by means of the Lie group action of
the structure group $GL(n,\mathbb R)$ on $\mathbb C$ given by
$(g,z)\mapsto |\det(g)|^\gamma z$. The smooth sections of
$|\!\bigwedge\!|^\gamma T^*Q$ are called complex $\gamma$-densities. 
In case $\gamma$ is real the bundle $|\!\bigwedge\!|^\gamma T^*Q$ contains
an obvious real subbundle (where the above action of the general linear
group is restricted to the real numbers in $\mathbb C$) which we shall call
the bundle of real $\gamma$-densities. The particular cases $\gamma=1$
and $\gamma=\frac{1}{2}$, respectively, are called the bundle of 
(real or complex) densities and half-densities, respectively. By a 
standard partition of unity argument there always exist nonvanishing 
sections of the bundle of real densities and hence the density bundles 
are trivial. Let us fix once and for all an arbitrary
nonvanishing real positive density $\mu$ on $Q$. This density 
locally defines a Riemann integral 
$\phi \mapsto \int \phi \mu_{1 \cdots n} dq^1 \cdots dq^n$
for any continuous real-valued function $\phi$ with support in a
coordinate neighbourhood $U$ in which $\mu$ takes the local form
$\mu_{1 \cdots n} dq^1 \cdots dq^n$, which is extended to all of $Q$
by means of the usual partition of unity argument and Riesz' theorem
to a Lebesgue integral on $Q$.

The connection $\nabla$ now defines a unique one-form $\alpha$ by
\BEQ {AlphaDef}
                 {\nabla}_X\mu =:\alpha(X)\mu
\EEQ
for an arbitrary vector field $X$ on $Q$ where the extension of $\nabla$
to the bundle of $\gamma$-densities is obvious. For any two vector fields
$X,Y$ on $Q$ we get upon evaluating 
$\nabla_X \nabla_Y\mu - \nabla_Y \nabla_X \mu - \nabla{[X,Y]} \mu $ 
the formula
\BEQ {CurvDAlpha}
          \tr R(X,Y)= -d\alpha(X,Y) .
\EEQ
Thus having fixed the one-form $\alpha$ we can now take the operator
$N$ as in (\ref{DefNeumaier}) and pass to the star product of Weyl 
type $\starw$ given by (\ref{DefStarw}).

Let $C^\infty_Q(T^*Q)[[\lambda]]$ denote the space of all formal series
in $\lambda$ whose coefficients lie in the space of all those 
smooth complex-valued functions on $T^*Q$ whose support intersected with
(the zero-section) $Q$ is compact. This is clearly a two-sided ideal of 
$(C^\infty(T^*Q)[[\lambda]], \starw)$ stable under complex conjugation. 
Moreover, let $C^\infty_0 (Q)[[\lambda]]$ denote the space of all 
formal series in $\lambda$ whose coefficients are smooth complex-valued 
functions of compact support. This space carries a 
$\mathbb C [[\lambda]]$-sesquilinear form defined by
\BEQ {DefSesqMu}
  \SP \phi \psi := \int_Q \cc{\phi} \psi \: \mu .
\EEQ
\begin{LEMMA} \label {SymStandWeylLem}
The formal series of differential operators corresponding to the 
standard and Weyl representation $\srep$ and $\wrep$ enjoy the following
symmetry properties for $f \in C^\infty(T^*Q)[[\lambda]]$ and
$\phi, \psi \in C^\infty_0 (Q)[[\lambda]]$:
\BEQ {SymStandWeyl}
     \SP {\srep (f) \phi} {\psi}
     = \SP {\phi} {\srep(N^2 \cc f)\psi}   
     \quad \mbox { and } \quad
     \SP{\wrep(f)\phi}{\psi}
     = \SP {\phi} {\wrep(\cc f)\psi}.      
\EEQ
\end{LEMMA}
\begin{PROOF}
The first equation is proved by integration by parts analogously to
the proof given in \cite[Theorem 7.3]{BNW97a} where this equation
had been proved for functions $\phi,\psi$ having their support in
a coordinate neighbourhood. By writing the coefficients of 
$\phi$ and $\psi$ in each order of $\lambda$ as finite
linear combinations of smooth complex-valued functions with support 
in coordinate neighbourhoods of an atlas and by (anti-) linearity the
local proof can be extended to the above statement. 
Then the second equation is an immediate consequence of the first and
the definition of $\wrep$.
\end{PROOF}

We are now ready to define and discuss a GNS construction analogous to
the one considered for flat $\mathbb R^{2n}$ in 
\cite[Proposition 11]{BW96b}. The positive functional which we 
want to consider is just the integration over the configuration space with 
respect to the fixed densitiy $\mu$:
\begin{PROPOSITION} \label{DefOmegaSatz}
Let $\omega_\mu : C^\infty_Q (T^*Q)[[\lambda]] \to \mathbb C[[\lambda]]$ 
be the following $\mathbb C[[\lambda]]$-linear functional
\BEQ {DefOmega}
    \omega_\mu (f):=\int_Q i^*f \: \mu 
\EEQ
then $\omega_\mu$ is well-defined and we have
\BEQ {NInv}
    \omega_\mu (N f) = \omega(f)  
\EEQ
\BEQ {NInvZwei}    
    \omega_\mu (f \starw g) 
    = \int_Q i^*(N^{-1}f)i^*(Ng) \: \mu .
\EEQ
It follows that $\omega_\mu$ is formally positive, i.~e.:
\BEQ {OmegaPositiv}
    \omega_\mu (\cc{f}\starw f)  = \int_Q i^*(\cc{Nf})i^*(Nf) \: \mu 
    \ge 0
\EEQ
and its Gel'fand ideal is given by
\BEQ {GelfandOmega}
    \mathcal J_\mu = \left\{ f \in C^\infty_Q(T^*Q)[[\lambda]] 
                     \; \big| \; i^*Nf = 0 \right\} .
\EEQ
\end{PROPOSITION}
\begin{PROOF}
To prove (\ref {NInv}), (\ref {NInvZwei}) and (\ref {OmegaPositiv})
it is sufficient to consider functions in $f, g \in C^\infty_Q (T^*Q)$ 
due to lemma \ref {ExtendPosFunctLem}.
Since $i^*f$ has compact support it follows that $\omega_\mu$ is
well-defined. Let $\chi$ be a smooth function with 
values in $[0, 1]$ and compact support on $Q$ which is equal to $1$ 
on an open neighbourhood containing the support of $i^*f$. Clearly, 
\[
    i^*f = \chi i^*f \chi = \chi i^*(f \pi^*\chi) = 
           \chi i^*(f \stars \pi^*\chi) = \chi \srep(f)(\chi)
\]
where the second to last equality follows from the fact that the vertical
derivatives of $f$ at the zero section vanish outside the intersection
of the support of $f$ with $Q$ whereas the covariant derivatives of
$\chi$ vanish on the intersection of the support of $f$ with $Q$.
Integrating the above equation over $Q$ with respect to $\mu$ 
and using lemma \ref {SymStandWeylLem} we get
\[
    \omega_\mu (f)  =  \SP{\chi} {\srep(f)\chi} 
                    =  \SP{\srep(N^2\cc{f})\chi}{\chi} 
                    =  \int_Q \srep(N^{-2}f)(\chi) \chi \: \mu 
                    =  \int_Q i^*(N^{-2}f) \: \mu 
                    =  \omega_\mu (N^{-2}f) .
\]
Hence (using the fact that $N$ is invertible) 
$\omega_Q((N-\id)(N+\id)f) = 0$ and since $N$ starts with $\id$ it 
follows that $\id+N$ is still invertible on 
$C^\infty_Q(T^*Q)[[\lambda]]$ whence $\omega_\mu ((N-\id) f)=0$ 
for all $f\in C^\infty_Q(T^*Q)[[\lambda]]$ which proves  
(\ref {NInv}). For equation (\ref {NInvZwei}) 
we get using the first one and the definition of $\starw$ that 
\begin {eqnarray*}
    \omega_\mu (f\starw g) 
    & = & \int_Q i^*((Nf) \stars (Ng)) \chi \chi \: \mu
          = \int_Q i^*( (Nf) \stars (Ng) \pi^* \chi) \chi \: \mu  \\
    & = & \int_Q i^*(((Nf) \stars (Ng)) \stars \pi^* \chi) \chi \: \mu
          = \SP{\chi} {\srep( (Nf)\stars (Ng) )\chi} \\
    & = & \SP {\srep(Nf)^\dag\chi}{\srep(Ng) \chi}
          = \SP {\srep (N\cc{f})\chi}{\srep(Ng) \chi} \\
    & = & \int_Q i^*(N^{-1}f)\chi i^*(Ng) \chi \: \mu
          = \int_Q i^*(N^{-1}f) i^*(Ng) \: \mu , 
\end{eqnarray*}
where $\chi$ was chosen as above obeying the condition to be equal to one
on an open set in $Q$ containing the union
of the supports of $i^*f$ and $i^*g$. Then equation (\ref{OmegaPositiv}) 
is an obvious particular case of equation (\ref{NInvZwei}) 
and implies that $\omega_\mu (\cc{f}\starw f) = 0$ iff $i^*N f=0$ 
vanishes since the integrand on the right hand side of 
equation (\ref{OmegaPositiv}) is formally positive: this proves the simple
formula for the Gel'fand ideal of $\omega_\mu$.
\end{PROOF}

Now we use the positive functional $\omega_\mu$ to define a GNS 
representation in the following standard way: Let 
$\mathfrak H_\mu := C^\infty_Q (T^*Q) [[\lambda]]/\mathcal J_\mu$ be the 
GNS (pre-) Hilbert space and denote the equivalence class
in $\mathfrak H_\omega$ of a function $f$ by $\psi_f$. Then the 
$\mathbb C[[\lambda]]$-valued Hermitian product in $\mathfrak H_\mu$ 
is given by $\SP {\psi_f} {\psi_g} := \omega_\mu (\cc f \starw g)$ 
and the representation $\pi_\mu (f)$ is determined by 
$\pi_\mu (f) \psi_g := \psi_{f \starw g}$. First we notice that the 
representation $\pi_\mu$ is not only defined for functions in 
$C^\infty_Q (T^*Q)[[\lambda]]$ but for all functions 
$C^\infty (T^*Q)[[\lambda]]$ since $\mathcal J_\mu$ is even a left 
ideal in $C^\infty (T^*Q)[[\lambda]]$ which can be verified directly
using the Cauchy-Schwartz inequality and the fact that 
$C^\infty_Q (T^*Q)[[\lambda]]$ is a two-sided ideal stable under 
complex conjugation. Now we can easily state the main theorem of 
this section characterising the representation more explicitly:
\begin {THEOREM} \label {GNSRepTheo}
With the notations from above we have:
\begin {enumerate}
\item The representation space $\mathfrak H_\mu$ is canonically isometric
      to $C^\infty_0 (Q)[[\lambda]]$ equipped with the 
      $\mathbb C[[\lambda]]$-valued Hermitian product (\ref{DefSesqMu}) 
      via the isometric isomorphism
      \BEQ {GNSHilbertIsom}
          \Phi: \psi_f \mapsto i^* Nf
          \quad \mbox { and its inverse } \quad
          \Phi^{-1}: \chi \mapsto \psi_{\pi^* \chi}
      \EEQ
      where $f \in C^\infty_Q (T^*Q)[[\lambda]]$ and 
      $\chi \in C^\infty_0 (Q)[[\lambda]]$.
\item The GNS representation $\pi_\mu$ carried to 
      $C^\infty_0 (Q)[[\lambda]]$ via $\Phi$ is just the Weyl 
      representation $\wrep$:
      \BEQ {GNSRepWrep}
          \Phi \pi_\mu (f) \Phi^{-1} \chi = i^* N (f \starw \pi^*\chi) 
          = \wrep (f) \chi. 
      \EEQ
\end {enumerate}
\end {THEOREM}       
\begin {PROOF}
Using the preceding proposition the proof is very simple: First one 
observes that $\Phi$ is indeed well-defined and bijective due to 
(\ref {GelfandOmega}) with inverse as stated above. Then 
\[
    \SP {\psi_f} {\psi_g} = \omega (\cc f \starw g) =
    \int_Q i^* \cc{(Nf)} i^* (Ng) \: \mu 
    = \int_Q \cc{\Phi (\psi_f)} \Phi (\psi_g) \: \mu = 
    \SP {\Phi(\psi_f)} {\Phi(\psi_g)}
\]
due to (\ref {NInvZwei}) which proves that $\Phi$ is isometric.    
Moreover we simply calculate 
\[
    \Phi \pi_\mu (f) \Phi^{-1} (\chi) = 
    \Phi \pi_\mu (f) \psi_{\pi^*\chi} =
    \Phi \psi_{f \starw \pi^*\chi} =
    i^* N( f \starw \pi^* \chi) 
\]
and thus (\ref {GNSRepWrep}) directly follows form the definition 
(\ref {wrepDef}) of $\wrep$.
\end {PROOF}

Notice that in view of appendix \ref {CNPApp} the whole 
construction can be done in the setting of formal Laurent or CNP series 
as well which justifies the terminology of `Hilbert spaces' and 
`Hermitian products' in sense of the definitions in 
\cite [Appendix A, B]{BW96b}, in particular $\CNP{C^\infty_0 (Q)}$ 
was shown to be already Cauchy-complete with respect to the topology
induced by $\SP \cdot \cdot$ in \cite [Theorem 7] {BW96b}.

\section {Several simple physical applications}
\label {ApplSec}

\subsection* {Convergence properties}

We shall briefly consider the question of convergence if one substitutes 
the formal parameter $\lambda$ by the numerical value 
$\hbar \in \mathbb R^+$ of Planck's constant (after having chosen some 
suitable unit system). In general this is a quite tricky problem and there 
has been a lot of work (and success!) to understand the convergence 
properties of the (a priori) formal star products, see e.~g.
\cite {CGRII,Fed96,BW96b,BBEW96b,Pfl97}. In the situation of homogeneous 
star products on a cotangent bundle the problem is luckily almost trivial 
since one has the $\mathbb C[\lambda]$-submodule 
$C^\infty_{pp} (T^*Q)[\lambda]$ of those functions which are polynomial in 
the momenta. Here the substitution $\lambda \mapsto \hbar$ makes no 
trouble at all since $\lambda$ appears only polynomially and furthermore 
the representations $\srep$ and $\wrep$ make $C^\infty_0 (Q)[\lambda]$ to 
a $C^\infty_{pp} (T^*Q)[\lambda]$-submodule of the 
$C^\infty (T^*Q)[[\lambda]]$-module $C^\infty_0 (Q)[[\lambda]]$
and thus again there is no problem to substitute $\lambda$ by $\hbar$ in 
the representations. We shall formulate this obvious fact not as a lemma 
but mention that the functions polynomial in the momenta are the most 
interesting observables for the physicist since in many examples the 
typical Hamiltonians and the typical integrals of motion are of this type.

\subsection* {Physical Examples}

Now we want to consider additional structures on $T^*Q$ which are 
motivated by various physical situations: Usually the Hamiltonian of a 
free particle moving on $Q$ is determined by a Riemannian metric $g$ on 
$Q$ by (setting the mass of the particle equal to $1$)
\BEQ {FreeParticle}
    \Hfree (q, p) := \frac{1}{2} g^\sharp (q) (p, p)
                         = \frac{1}{2} g^{ij}(q) p_i p_j
\EEQ
where $g^\sharp \in \Gamma (\bigvee^2 TQ)$ is the `inverse' metric tensor.
Obviously $\Hfree \in C^\infty_{pp,2} (T^*Q)$ where the corresponding 
tensor field is just $g^\sharp \in \Gamma (\bigvee^2 TQ)$.
In this case we shall use the Levi-Civita connection 
$\nabla^{\mbox{\tiny LC}}$  of $g$ to construct the star product $\stars$ 
and since the metric induces a canonical volume density $\mu_g$ we shall 
use this density to define the one-form $\alpha$, the corresponding 
operator $N$, and the star product of Weyl ordered type $\starw$. Now 
$\nabla^{\mbox{\tiny LC}} \mu_g = 0$ and hence $\alpha = 0$ which leads to 
\BEQ {NeumaierMetric}
    N = e^{\frac{\lambda} {2\im} \Delta}
    \qquad
    \mbox { with }
    \qquad
    \Delta = \partial_{p_i} \partial_{q^i} 
           + p_r (\pi^* \Gamma^r_{ij}) \partial_{p_i} \partial_{p_j}
           + (\pi^*\Gamma^i_{ij}) \partial_{p_j} .           
\EEQ
Furthermore it is easy to calculate using (\ref {wrepExpl}) that the 
operator corresponding to $\Hfree$ is given by
\BEQ {WeylRepHfree}
    \wrep (\Hfree) = - \frac{\lambda^2} {2} \Delta_g
\EEQ
where $\Delta_g$ is the Laplacian of the metric $g$. 
Moreover {\em no} `quantum potentials' in $\wrep (\Hfree)$ 
like e.~g. a multiple of the Ricci scalar occur as in the approach 
of Underhill \cite {Und78}. If in addition a
smooth potential $V \in C^\infty (Q)$ is present then the physical 
Hamiltonian will be given by
\BEQ {HfreeV}
    H = \Hfree + \pi^* V
\EEQ
and the corresponding operator is 
\BEQ {WeylRepHfreeV}        
    \wrep (H) = -\frac{\lambda^2}{2} \Delta_g + V
\EEQ
where $V$ acts simply as left multiplication. Another important 
physical example is given by functions linear in the momenta since 
they generate the point transformations of the configuration space $Q$. 
Let $\widehat{X} \in C^\infty_{pp,1} (T^*Q)$ be a function linear in the 
momenta and let $X \in \Gamma (TQ)$ be the corresponding vector field then
we obtain
\BEQ {WeylRepVector}
    \wrep (\widehat{X}) = \frac{\lambda}{\im} \left(
                          \Lie_X + \frac{1}{2} \Div_g X \right)
\EEQ
where $\Lie_X$ denotes the Lie derivative with respect to $X$ and 
$\Div_g X$ denotes the metric divergence of $X$. 
Note that any homogeneous star product of Weyl type $*$ is 
covariant
under $C^\infty_{pp, 1} (T^*Q)$, i.~e.
$\widehat X * \widehat Y - \widehat Y * \widehat X 
= \im\lambda \{\widehat X, \widehat Y \}$ for all 
$X, Y \in \Gamma (TQ)$ (see e.~g. \cite {ACMP83} for further definitions 
of different types of invariance).

Summarizing we observe that the well-known `ad-hoc quantization rules' 
on cotangent bundles are obtained by deformation quantization and 
GNS representation in a very systematic way.

\subsection* {Symmetries of $T^*Q$}

At last we shall discuss the action of geometrical symmetries of $Q$ and 
$T^*Q$ as (anti-) automorphisms of the star product algebra 
$C^\infty (T^*Q) [[\lambda]]$. Firstly we shall remember some easy facts 
about the way (anti-) automorphisms induce (anti-) unitary maps via the 
GNS construction in a more general setting adapted form the usual 
$C^*$-theory (see appendix \ref {CNPApp} and \cite {BW96b} for 
definitions, notation and further properties).
\begin {PROPOSITION} \label {AutoGNSProp}
Let $\mathcal A, \mathcal B$ be associative $^*$-algebras over 
$\field C := \field R(\im)$ where $R$ is an ordered field and 
$\im^2 = -1$ and $^*$ is an involutive $\field C$-antilinear 
(with respect to the complex conjugation in $\field C$) 
anti-automorphism of $\mathcal A$ resp. $\mathcal B$. Moreover let 
$\omega: \mathcal B \to \field C$ be a positive linear functional with 
Gel'fand ideal $\mathcal J_\omega$ and GNS representation $\pi_\omega$ on 
$\mathfrak H_\omega := \mathcal B / \mathcal J_\omega$. Furthermore let 
$A: \mathcal A \to \mathcal B$ be a $^*$-homomorphism and 
$\tilde A: \mathcal A \to \mathcal B$ be a $^*$-anti-homomorphism.
Then we have:
\begin {enumerate}
\item The functionals $\omega_A := \omega \circ A$ and 
      $\omega_{\tilde A} := \omega \circ \tilde A$ are positive linear 
      functionals of $\mathcal A$ with Gel'fand ideals 
      $\mathcal J_A = A^{-1} (\mathcal J_\omega)$ and
      $\mathcal J_{\tilde A} = \tilde A^{-1} ({\mathcal  J_\omega}^*)$.
      
\item Let $\pi_A$ resp. $\pi_{\tilde A}$ be the 
      GNS representations on $\mathfrak H_A := \mathcal A / \mathcal J_A$
      resp. $\mathfrak H_{\tilde A} := \mathcal A / \mathcal J_{\tilde A}$
      induced by $\omega_A$ resp. $\omega_{\tilde A}$ then the maps
      \BEQ {AntiUnitaryDef}
          U_A : \begin {array} {c}
                \mathfrak H_A \to \mathfrak H_\omega \\
                \psi^A_f \mapsto \psi_{Af}
                \end {array}
          \quad
          \mbox { resp. }
          \quad
          U_{\tilde A} : 
                \begin {array} {c}
                \mathfrak H_{\tilde A} \to \mathfrak H_\omega \\
                \psi^{\tilde A}_f \mapsto \psi_{\tilde A (f^*)}
                \end {array}
      \EEQ
      (where $\psi_\cdot$, $\psi^A_\cdot$, and $\psi^{\tilde A}_\cdot$
      denote the equivalence classes in $\mathfrak H_\omega$, 
      $\mathfrak H_A$, and $\mathfrak H_{\tilde A}$)
      are well-defined and isometric resp. anti-isometric and one has
      \BEQ {UAGNSRep}
          U_A \pi_A (f) = \pi_\omega (Af) U_A
          \quad 
          \mbox { resp. }
          \quad
          U_{\tilde A} \pi_{\tilde A} (f) 
          = \pi_\omega (\tilde A f^*) U_{\tilde A} .
      \EEQ          
\item If moreover $A$ resp. $\tilde A$ is surjective then $U_A$ resp. 
      $U_{\tilde A}$ is unitary resp. anti-unitary and the inverse of $U_A$ 
      resp. $U_{\tilde A}$ is (well-defined!) given by
      \BEQ {InverseUA}
          U_A^{-1} \psi_g = \psi^A_f
          \quad
          \mbox { resp. }
          \quad
          U_{\tilde A}^{-1} \psi_g = \psi^{\tilde A}_f
      \EEQ
      where $f \in A^{-1} (\{g\})$ resp. $f \in {\tilde A}^{-1} (\{g^*\})$.         
\end {enumerate}
\end {PROPOSITION}
\begin {PROOF}
The proof is a straight forward computation using the (anti-) homomorphism 
property of $A$ (resp. $\tilde A$) as well as the compatibility of $A$ and 
$\tilde A$ with the $^*$-involution. 
\end {PROOF}

A homomorphism which respects the involution will also be called a real 
homomorphism.
\begin {COROLLARY} \label {InducedUnitaryCor}
With the notations from above let 
$A$ (resp. $\tilde A) : \mathcal A \to \mathcal A$ be a real 
(anti-) automorphism 
and let $\omega: \mathcal A \to \field C$ be a positive $A$-invariant
(resp. $\tilde A$-invariant) linear functional, i.~e. 
$\omega = \omega \circ A$ (resp $\omega = \omega \circ \tilde A$). 
Then $U_A$ 
(resp. $U_{\tilde A}) : \mathfrak H_\omega \to \mathfrak H_\omega$ is a 
(anti-) unitary map with inverse $U_A^{-1} \psi_f = \psi_{A^{-1} f}$ resp.
$U_{\tilde A}^{-1} \psi_f = \psi_{{\tilde A}^{-1} f^*}$ and
\BEQ {UARepUA}
    \pi_\omega (f) = U_A^{-1} \pi_\omega (Af) U_A
    \quad
    \mbox { resp. }
    \quad
    \pi_\omega (f) = 
    U_{\tilde A}^{-1} \pi_\omega (\tilde A f^*) U_{\tilde A} .
\EEQ    
\end {COROLLARY}
Note that all these results are still correct if one considers an 
ordered ring $\field R$ instead of an ordered field as we shall do in the 
following using $\mathbb R[[\lambda]]$ and $\mathbb C[[\lambda]]$.

Now we come back to the particular situation $T^*Q$. Let $\phi: Q \to Q$ 
be a diffeomorphism of $Q$ leaving invariant the connection $\nabla$, 
i.~e. $\nabla_X Y = \phi^* (\nabla_{\phi_* X} \phi_* Y)$ for all vector 
fields $X, Y \in \Gamma (TQ)$. Moreover we consider the canonical lift 
$T^*\phi : T^*Q \to T^*Q$ of $\phi$ to a symplectic diffeomorphism of 
$T^*Q$ (where we use the convention such that 
$T^*\phi \circ i = i \circ \phi$ and not $i \circ \phi^{-1}$).
\begin {LEMMA}
Let $\nabla$ be a $\phi$-invariant torsion-free connection on $Q$ where 
$\phi: Q \to Q$ is a diffeomorphism then the pull back 
$A_\phi := (T^*\phi)^* : 
C^\infty (T^*Q)[[\lambda]] \to C^\infty (T^*Q)[[\lambda]]$
is a real automorphism of the corresponding star product $\stars$, i.~e.:
\BEQ {AphiAutoStars}
    A_\phi (f \stars g) = A_\phi f \stars A_\phi g .
\EEQ
\end {LEMMA}    
\begin {PROOF}
We shall only indicate the proof very briefly: First we notice that in the 
Fedosov setting as desribed in \cite {BNW97a} the map $A_\phi$ naturally 
extends to the whole Fedosov algebra and defines a real automorphism of 
the fibrewise standard ordered product. This is a straight forward 
computation using the invariance of $\nabla$ under $\phi$. Next one proves 
that $A_\phi$ commutes with the Fedosov derivation $\Dstd$ using the 
recursion formulas for the curvature part of $\Dstd$. Finally one applies 
proposition 2.3 in \cite {BNW97a} and the rather obvious fact that 
$A_\phi$ commutes with the Fedosov Taylor series to prove the fact that 
the fibrewise automorphism induces indeed an automorphism of the star 
product, namely $A_\phi$. Another possibility is given by proving 
(\ref {AphiAutoStars}) by first restricting to functions polynomial 
in the momenta which turns out to be sufficient to prove 
(\ref {AphiAutoStars}) and secondly using the explicit formula for 
the representation $\srep$ and the invariance of $\nabla$.
\end {PROOF}

Now if $A_\phi$ is an automorphism of $\stars$ it is an automorphism of 
$\starw$ too, if it commutes with the operator $N$ which is clearly the 
case if the one-form $\alpha$ is $\phi$-invariant: 
$\phi^* \alpha = \alpha$. Note that a priori we only know from the 
$\phi$-invariance of $\nabla$ that $d\alpha = - \tr R$ is 
$\phi$-invariant. But assuming $\phi^*\alpha = \alpha$ we obtain that 
$A_\phi$ is an automorphism of $\starw$ as well. Note that in the case 
where $\nabla$ is uni-modular, i.~e. $\tr R = 0$ (and in particular in 
the Riemannian case with the Levi-Civita connection), 
with the choice $\alpha = 0$ this will automatically be fulfilled.
\begin {COROLLARY}
Let $\alpha \in \Gamma (T^*Q)$ be a one-form such that $d\alpha = -\tr R$ 
and $\phi^*\alpha = \alpha$ then $A_\phi = (T^*\phi)^*$ is a real 
automorphism of the corresponding star product $\starw$ and commutes 
with $N$.
\end {COROLLARY}
Assuming in addition that we have a $\phi$-invariant volume density $\mu$ 
we automatically obtain that $\alpha$ defined by (\ref {AlphaDef}) is 
$\phi$-invariant. Moreover $\omega_\mu$ defined as in (\ref {DefOmega}) is 
an $A_\phi$-invariant positive linear functional and hence we can apply 
corollary \ref {InducedUnitaryCor} to obtain the following easy lemma:
\begin {LEMMA}
Let both $\nabla$ and $\mu$ be $\phi$-invariant then 
$\phi^* \alpha = \alpha$ and $A_\phi$ is a real $\starw$-automorphism. 
Furthermore $\omega_\mu$ is $A_\phi$-invariant and the induced unitary map 
$U_\phi$ in the GNS Hilbert space is simply given by
\BEQ {SymmUnitary}
    U_\phi \chi = \phi^* \chi 
    \qquad
    \chi \in C^\infty_0 (Q)[[\lambda]] .
\EEQ    
\end {LEMMA}
This lemma covers obviously all kinds of (Lie-) group actions on $Q$ which 
leave invariant a connection and a volume density and ensures that such 
classical group actions are implemented as group actions of automorphisms 
of the observable algebra as well as unitray group actions on the GNS 
Hilbert space both realized by pull backs. This is of course a well-known 
implementation but again deformation quantization and the 
GNS representation provide a very systematical way to study such 
geometric symmetries and their realizations in quantum mechanics.

At last we shall discuss the `time reversal' on a kinematical level, i.~e. 
as a geometrical property of the classical phase space 
(see e.~g. \cite [p. 308] {AM85}):
\begin {DEFINITION}
Let $(M, \omega)$ be a symplectic manifold and $T: M \to M$ a 
diffeomorphism then $T$ is called time reversal map if 
$T^*\omega = - \omega$.
\end {DEFINITION}
\begin {LEMMA}
If $(M,\omega)$ has a time reversal map $T$ then $T \circ \phi$ is again a 
time reversal map for any symplectic diffeomorphism $\phi$ and any time 
reversal map of $M$ is obtained this way. If $\omega$ is not exact then 
$T$ is a `large diffeomorphism' i.~e. not isotopic to the identity. 
If the fixed point set of $T$ is a submanifold then it is isotropic.
\end {LEMMA}
In the case of a cotangent bundle we canonically have such a time reversal 
map namely $T: T^*_q Q \ni \alpha_q \mapsto -\alpha_q$ with the additional
property $T^2 = \id$. In a local bundle Darboux chart this reads 
$T: (q, p) \mapsto (q, -p)$ which justifies the name `time reversal'. 
The following proposition should be well-known:
\begin {PROPOSITION}
Let $(M, \omega)$ be a symplectic manifold and $T: M \to M$ a time 
reversal map with $T^2 = \id_M$. Then the set of fixed points of $T$ is 
either empty or a (not necessarily connected) Lagrangean submanifold of 
$M$.
\end {PROPOSITION}
\begin {PROOF}
First we notice that there always exists a symplectic torsion-free and 
$T$-invariant connection $\nabla$. Now assume that 
$L := \{ p \in M | T(p) = p\}$ is not empty and let $p \in L$. Then 
$T_p T: T_pM \to T_p M$ has square $\id_{T_pM}$ and hence it is 
diagonalizable with eigenspaces $E_p^\pm$ to the eigenvalues $\pm 1$. 
Clearly $E_p^{\pm}$ are both Lagrangean subvector spaces of $T_p M$. 
Now the exponential mapping $\exp_p$ of $\nabla$ is locally a 
diffeomorphism and maps a neigbourhood of $0_p$ in $E_q^+$ into $L$ due 
to the $T$-invariance of $\nabla$ and thus the local inverse of $\exp_p$ 
is a submanifold chart for $L$ in a 
neighbourhood of $p$ proving thereby that $L$ is indeed a submanifold 
(even totally geodesic with respect to $\nabla$) with tangent space 
$E_p$ at $p$ for all $p \in L$. This implies that $L$ is Lagrangean.
\end {PROOF}

Certain coadjoint orbits carry time reversal maps with nonempty
fixed-point set:
\begin {PROPOSITION}
Let $\mathfrak g$ a real finite-dimensional Lie
algebra and $\mathfrak g^*$ its dual space. For any
$\mu\in \mathfrak g^*$ consider the coadjoint orbit $\mathcal O_\mu$
through $\mu$, i.~e.
$\mathcal O_\mu := \{ \exp(\ad^* (\xi_1)) \cdots \exp(\ad^* (\xi_n)) \mu
\in \mathfrak g^* \; | \;  n \in \mathbb N, \xi_1, \ldots, \xi_n \in
\mathfrak g\}$.
\begin{enumerate}
\item Let furthermore
      $s: \mathfrak g \to \mathfrak g$ an involutive
      automorphism of Lie algebras (i.~e. $s^2 = \id_{\mathfrak g}$).
      Then $\mathfrak g^*$ decomposes into the direct sum
      of the eigenspaces $\mathfrak k^*$ and
      $\mathfrak m^*$ of $s^*$, the dual map to $s$,
      corresponding to the eigenvalues $1$ and $-1$, respectively.
      Pick $\mu \in \mathfrak m^*$ and consider the coadjoint orbit
      $\mathcal O_\mu$ through $\mu$.
      Then the map $T: \mathcal O_\mu \to \mathcal O_\mu$
      defined by the restriction of $-s^*$ is an involutive
      time reversal map of the symplectic manifold $\mathcal O_\mu$ with
      its Kirillov-Kostant-Souriau form whose fixed point set is equal to
      the intersection of $\mathfrak m^*$ and $\mathcal O_\mu$.

\item For each element $\mu$ in the dual space $\mathfrak g^*$ of a
      semisimple compact Lie algebra $\mathfrak g$
      there is an involutive automorphism $s$ of $\mathfrak g$ such
      that $\mu$ is in the eigenspace $\mathfrak m^*$ of $s^*$.
      In view of {\it i.)} it follows that every coadjoint orbit
      in $\mathfrak g^*$ with semisimple compact $\mathfrak g$
      admits an involutive time reversal map with nonempty fixed point
      set.
\end{enumerate}
\end{PROPOSITION}
\begin {PROOF}
  {\it i.)} Since $s^* \ad^* (\xi)=ad^* (s \xi) s^*$
  it follows that $-s^*$
  restricts to $\mathcal O_\mu$ whence $T$ is well-defined and the
  fixed point set is equal to the intersection of the orbit with
  $\mathfrak m^*$. Moreover, this implies that for all
  $\xi, \eta \in \mathfrak g$ and $\nu \in \mathcal O_\mu$ we have
  $T_\nu T \xi_{\mathcal O_\mu} (\nu) = (s \xi)_{\mathcal O_\mu} (T \nu)$
  for all the vector fields $\xi_{\mathcal O_\mu}(\nu) = \ad^* (\xi) \nu$
  implying that
  $(T^* \omega)(\nu)(\xi_{\mathcal O_\mu}(\nu),
  \eta_{\mathcal O_\mu}(\nu)) = - (s^* \nu)([s \xi, s \eta])
  = - \nu([\xi, \eta])
  = - \omega(\nu)(\xi_{\mathcal O_\mu}(\nu), \eta_{\mathcal O_\mu}(\nu))$.

  {\it ii.)} Identifying the Lie algebra and its dual space by means of
  the Killing form we can put $\mu \in \mathfrak g$ in some
  Cartan subalgebra $\im \mathfrak h$ of $\mathfrak g$. Using the
  root space decomposition
  \[
   \mathfrak g = \sum_{\alpha \in \Delta}
                 \mathbb R (X_\alpha - X_{- \alpha})~
                 \oplus~\im\mathfrak h~
                 \oplus~\sum_{\alpha \in \Delta}
                 \mathbb R\im (X_\alpha + X_{-\alpha})
  \]
  where $\Delta$ is the set of all roots of $\mathfrak g$, and $X_\alpha$
  are normalized vectors in the eigenspaces $\mathfrak g_\alpha$
  (see \cite [p.~182] {Hel78} for definitions and details)
  we see that the first summand is a subalgebra $\mathfrak k$ of
  $\mathfrak g$ whereas the second plus the third summand is a subspace
  $\mathfrak m$ of $\mathfrak g$ such that
  $[\mathfrak k, \mathfrak k] \subseteq \mathfrak k$,
  $[\mathfrak k, \mathfrak m] \subseteq \mathfrak m$, and
  $[\mathfrak m, \mathfrak m] \subseteq \mathfrak k$, where we have
  made use of \cite [Theorem 5.5]{Hel78}. Hence $\mathfrak g$ is an
  orthogonal symmetric Lie algebra
  (see \cite [p.~213] {Hel78}) where $s$ can be defined by
  being $1$ on $\mathfrak k$ and $-1$ on $\mathfrak m$.
\end{PROOF}

Now we come back to cotangent bundles and consider a Fedosov star product 
of Weyl ordered type $\starw$ on $T^*Q$ then we have the following lemma 
which is proved directly using the Weyl type property and the homogeneity 
of $\starw$:
\begin {LEMMA}
The canonical time reversal map $T$ induces via pull back a real 
anti-automorphism of any $\starw$ on $T^*Q$, i.~e. for any 
$f, g \in C^\infty (T^*Q)[[\lambda]]$ and $A_T := T^*$ we have
\BEQ {TimeRevAntiAuto}
    A_T (f \starw g) = A_T g \starw A_T f .
\EEQ
\end {LEMMA}
Since $A_T$ is an anti-automorphism of $\starw$ we obtain an anti-unitary 
map $U_T$ in each GNS Hilbert space $\mathfrak H_\mu$ since obviously 
$\omega_\mu \circ A_T = \omega_\mu$ for all volume densities $\mu$. 
As expected it turns out that $U_T$ is just the complex conjugation 
of the wave functions:
\begin {LEMMA}
Let $A_T$ be the time reversal $\starw$-anti-automorphism and let 
$\omega_\mu$ be the positive linear functional as in 
(\ref {DefOmega}) then the induced anti-unitary map $U_T$ in the GNS 
Hilbert space $C^\infty_0 (Q)[[\lambda]]$ is given by complex conjugation
\BEQ {TimeRevCC}
    U_T \chi = \cc \chi 
    \qquad
    \chi \in C^\infty_0 (Q)[[\lambda]] .
\EEQ
\end {LEMMA}

\section {Comparison with other approaches and integral formulas}
\label {CompSec}

In this section we shall compare our differential operator
representations we constructed to some different approaches using
integral representations. Before we can translate our explicit 
formulas to integral expressions we have to substitute the formal 
parameter $\lambda$ by the real positive number $\hbar \in \mathbb R^+$
which is done most easily by restricting to functions polynomial in the 
momenta. The substitution $\lambda \mapsto \hbar \in \mathbb R^+$ is 
denoted by $\ldots |_{\lambda = \hbar}$. In general the integral 
expressions will be well-defined for a larger class of functions 
namely certain symbol classes which are smooth functions on $T^*Q$ 
with controllable fibrewise increase. Then our formal results can be 
obtained by asymptotic expansions.

First of all we shall show that the
integral expression given by \cite{Pfl97} in the framework of 
symbol calculus when applied to functions that are polynomial in
the momenta coincides with the representation $\srep$ we
introduced in a purely algebraic fashion. Moreover we define some
integral representations that are quite natural generalizations of
the  Weyl quantization in the case of flat $T^*\mathbb R^n$
(c.~f. \cite{Fed96}) and have been already studied by \cite{Emm93a,Emm93b}
in the case of a Riemannian manifold $Q$.
In so far our results are generalizing since we drop those technically 
simplifying preconditions. On the other hand we can show that this 
alternative approach to a Weyl quantization directly yields (when 
applied to polynomial functions in the momenta) the Weyl representation 
we used to define our star product $\starw$ in order to do the GNS 
construction.

\subsection{An integral expression for the standard representation}

In the flat case it is well-known that the standard representation
$\srep(\widehat{T})$ of $\widehat{T}\in C_{pp}^\infty(T^*\mathbb R^n)$
applied to a function $\phi \in C_0^\infty(\mathbb R^n)$ can be given 
by the integral representation:
\[
    \left.(\srep(\widehat{T})\phi) (q)\right|_{\lambda = \hbar} 
    = \frac{1}{(2\pi \hbar)^n} \int_{\mathbb R^n}
      \widehat{T}(q,p) \int_{\mathbb R^n}
      e^{-\frac{\im}{\hbar}\langle p,v \rangle}\phi(q + v ) 
      d^n v d^n p,
\]
and therefore we first report on the generalization of such an
expression as it has been proposed in \cite{Pfl97} and give an
analogue that can be defined pointwise for a fixed $q \in Q$ by
integrations over $T_q^*Q$ and $T_qQ$ using the canonical
symplectic volume form $\Omega_q$ on $T_qQ \times T_q^*Q$. The
expression $\phi (q+ v)$ in the above formula obviously makes no
sense in the case of an arbitrary manifold $Q$ but can be viewed as
$\phi( \exp_q(v_q))$ for $v_q \in T_q \mathbb R^n$ denoting by
$\exp$ the exponential mapping with respect to the flat connection.
In general the exponential mapping not being globally defined we
aim to define such an expression just in a neighbourhood of $0_q$.
To do so, we choose an open neighbourhood $\mathcal O \subseteq TQ$
of the zero section that is mapped diffeomorphically to an open
neighbourhood of the diagonal in $Q \times Q$ by $\pi \times \exp$.
By possibly shrinking this neighbourhood we may assume
that for every $v_q \in \pi^{-1}(q) \cap \mathcal{O}$ 
also $- v_q$ is contained in $\pi^{-1}(q) \cap \mathcal{O}$. 
This will be important for the definition of the Weyl 
representation later. By Urysohn's lemma there is a smooth function 
$\chi: TQ \to [0,1]$ and an open neighbourhood 
$\widetilde{\mathcal{O}} \subseteq \mathcal{O}$
of the zero section such that $\chi|_{\widetilde{\mathcal{O}}}=1$ and
$\mathrm{supp} (\chi) \subseteq \mathcal{O}$. Using this so-called
cut-off function we may assign to any smooth function $\phi$ on $Q$
a function $\phi_q^\chi \in C_0^\infty(T_qQ)$ by
\[
    \phi_q^\chi(v_q) := \left\{ \begin{array} {rcl} 
                        \chi(v_q) \phi(\exp_q(v_q)) & 
                        \textrm{for} & 
                        v_q \in \pi^{-1}(q)\cap \mathcal{O}\\ 
                        0 & 
                        &
                        \textrm{else} .
                        \end{array} \right.
\]
After these preparations the mapping
\BEQ {standint}
    \mathcal S(\widehat{T}) : C^\infty(Q) \ni \phi 
    \mapsto \left( q \mapsto \frac{1}{(2\pi \hbar)^n} \int_{T^*_q Q} 
    \widehat{T}(\alpha_q)\int_{T_qQ} e^{-\frac{\im}{\hbar}\alpha_q(v_q)} 
    \phi_q^\chi(v_q) \Omega_q\right)
\EEQ
is obviously well-defined for $\widehat{T} \in C_{pp}^\infty(T^*Q)$ by
well-known properties of the Fourier transform. At first sight it
might seem that this expression depends on the special choice of
the function $\chi$, what actually is not the case for functions 
polynomial in the momenta, and in addition
it is a priori not clear whether the function 
$\mathcal S (\widehat{T})\phi$ is again a smooth function. 
But we have the following:
\begin{LEMMA}
For all $\widehat{T} \in C_{pp}^\infty(T^*Q)$ and $\phi \in
C^\infty(Q)$ we have
\[
\mathcal S(\widehat{T}) \phi = \srep(\widehat{T})\phi|_{\lambda=\hbar}.
\]
\end{LEMMA}
\begin{PROOF}
First of all we notice that $(\mathcal S(\widehat{T})\phi)(q)$ is 
defined independently on the choice of the basis of $T_qQ$ and $T^*_qQ$ one 
introduces to carry out the integrations, therefore it is convenient to use 
a coordinate basis induced by a normal chart around $q$.
In addition it suffices to prove the assertion for homogeneous 
$\widehat{T}$. By an obvious calculation one gets 
$(\mathcal S(\widehat{T})\phi)(q)= 
\left(\frac{\hbar}{\im}\right)^k \frac{1}{k!}T^{i_1 \ldots i_k}(q) 
\left.\frac{\partial^k \phi_q^\chi(v_q)} {\partial v^{i_1} \cdots 
\partial v^{i_k}}\right|_{v_q = 0_q}$. 
Now by construction of the 
cut-off function $\chi$ all its derivatives at $0_q$ 
vanish and we get the coordinate expression for 
$(\srep(\widehat{T})\phi)(q)|_{\lambda=\hbar}$ in a normal chart around $q$ 
proving the assertion since in a normal chart around $q$ the 
symmetric covariant derivative is given by 
$(D^k \phi)(q) = \frac {\partial^k \phi}
                 {\partial q^{i_1} \cdots \partial q^{i_k}}(q) dq^{i_1}
                 \vee \ldots \vee dq^{i_k}$ 
(c.~f. \cite[Lemma A.9]{BNW97a}).
\end{PROOF}

\subsection{An integral expression for the Weyl representation}

Again motivated by the flat case in \cite{Emm93a,Emm93b} an integral 
expression has been introduced to generalize the Weyl quantization in 
the following manner. For $\phi, \psi \in C_0^\infty(Q)$ and 
$\widehat{T}\in C_{pp}^\infty(T^*Q)$ one considers using the same notation 
as in the preceding subsection the mapping
\[
    \mathsf W(\widehat{T}) : (\phi, \psi) 
    \mapsto \left(q \mapsto \frac{1}{(\pi 
    \hbar)^n} \int_{T^*_qQ} \widehat{T}(\alpha_q) \int_{T_qQ} 
    e^{-\frac{2\im}{\hbar}\alpha_q(v_q)} \overline{\phi_q^\chi (-v_q)} 
    \psi_q^\chi(v_q) \Omega_q\right), 
\]
that shall serve as an integral kernel to define the Weyl representation 
$\mathcal W(\widehat{T}) \psi$ by its `Matrix elements' 
$\langle \phi, \mathcal W (\widehat{T})\psi\rangle := 
\int_Q  \mathsf W(\widehat{T}) (\phi, \psi)\,\mu$. 
This is clearly well-defined, since 
$(C_0^\infty (Q), \langle{}\cdot {},{}\cdot {}\rangle)$ is a 
pre-Hilbert space and as we shall see in the next lemma 
$\mathsf W(\widehat{T})$ is a differential operator acting on the 
functions $\phi,\psi$.
\begin{LEMMA}
$\mathsf W(\widehat{T})(\phi,\psi)$ is defined independently on the choice of 
the cut-off function $\chi$ and in a local chart of $Q$ it is given for 
$\widehat{T} \in C_{pp,k}^\infty(T^*Q)$ by
\begin{eqnarray*}
    \mathsf W(\widehat{T})(\phi,\psi)
    & = & \frac{1}{k!}
          \left(\frac{\hbar}{2\im} \right)^k 
          i_s^*(dq^{i_1}) \ldots i_s^*(dq^{i_k})T \\
    &   & \times \sum_{r=0}^k (-1)^r
          {k \choose r}
          i_s(\partial_{q^{i_1}}) \ldots i_s(\partial_{q^{i_r}}) 
          \frac{1}{r!} D^r \overline{\phi} \, 
          i_s(\partial_{q^{i_{r+1}}}) \ldots i_s(\partial_{q^{i_k}}) 
          \frac{1}{(k-r)!} D^{k-r} \psi,
\end{eqnarray*}
hence it obviously is smooth, and since the functions $\phi, \psi$ have  
compact support this implies that 
$\int_Q \mathsf W(\widehat{T})(\phi,\psi) \mu$ is well-defined.
\end{LEMMA}
\begin{PROOF}
By a straight forward computation as in the case of the standard 
representation using a normal chart around $q$ one obtains
$\mathsf W(\widehat{T})(\phi,\psi)(q) = 
\frac{1}{k!}\left(\frac{\hbar}{2\im}
\right)^k T^{i_1\ldots i_k}(q)\left.\frac{\partial^k 
(\overline{\phi_q^\chi(-v_q)} 
\psi_q^\chi (v_q))}{\partial v^{i_1} \cdots \partial 
v^{i_k}}\right|_{v_q=0_q}$. 
By the Leibniz rule and the construction 
of $\chi$ one gets the assertion observing the shape of the symmetric 
covariant derivative in a normal chart around $q$.
\end{PROOF}\\
After this preparation we can state the following Lemma:
\begin{LEMMA}
The Weyl representation $\mathcal W(\widehat{T})\psi$ being well-defined by 
the equation $\langle \phi, \mathcal W(\widehat{T})\psi \rangle = \int_Q 
\mathsf W(\widehat{T})(\phi,\psi) \,\mu$ coincides with the 
Weyl-representation $\wrep$ defined as in (\ref {wrepDef}) after the 
substitution $\lambda \mapsto \hbar$, i.~e.
\[
    \mathcal W(\widehat{T})\psi = 
    \left. \wrep(\widehat{T})\psi\right|_{\lambda=\hbar} = 
    \left. \srep(N \widehat{T})\psi\right|_{\lambda=\hbar},
\]
where $N$ is given as in (\ref{DefNeumaier}).
\end{LEMMA}
\begin{PROOF}
The proof is an easy but lengthy calculation using iterated 
integrations by parts in the same fashion as they were nessecary to prove 
the equations (\ref{SymStandWeyl}). 
\end{PROOF}\\

Finally, note that other approaches like Underhill's (see \cite{Und78}) 
use compactly 
supported half-densities as representation space instead of compactly 
supported functions on $Q$. But although there is an isomorphism
between these spaces by choosing a fixed positive density $\mu$ 
and assigning to $\phi \in C_0^\infty(Q)$ the half-density 
$\check{\phi}= \phi \mu^{1/2}$ the corresponding operators turn out to 
differ
e.~g. by additional terms proportional to $\hbar^2$ times the scalar
curvature for the free particl Hamiltonian(\ref{FreeParticle}): this 
difference is due to whether the
reference density $\mu^{1/2}$ is pulled-back by means of the
exponential map or not, see e.~g.~\cite[518--520]{Emm93a} for a detailed 
discussion.

\section {WKB expansion for projectable Lagrangean submanifolds}
\label {WKBSec}

In this section we shall discuss how the usual WKB expansion can be 
formulated in the framework of deformation quantization using particular
GNS representations in the case of cotangent bundles. We consider here 
a Hamiltonian $H \in C^\infty (T^*Q)$ and assume that for a fixed energy 
value $E \in \mathbb R$ in the image of $H$ there exists a Lagrangean 
submanifold $L_\beta$ contained in $H^{-1} (\{E\})$, i. e.
\BEQ {HamiltonJacobiEq}
    H | L_{\beta} = E
\EEQ
which is furthermore given by the graph of a closed one-form 
$\beta \in \Gamma (T^*Q)$, i. e.
\BEQ {LbetaDef}
    L_\beta = \graph (\beta) 
    = \left\{ \alpha_q \in T^*Q \; \big| \; \forall q \in Q: 
    \alpha_q = \beta (q) \right\} .
\EEQ
Such Lagrangean submanifolds are called projectable and it is well-known 
that $d\beta = 0$ is equivalent to the statement that $\graph (\beta)$ is 
Lagrangean, see e.~g. \cite [Prop. 5.~3.~15] {AM85}. For later use we 
denote by 
\BEQ {ibetaDef}
    i_\beta: L_\beta \to T^*Q
\EEQ
the embedding of $L_\beta$ in $T^*Q$. Since $d\beta = 0$ we have 
$d\pi^*\beta = 0$ which implies that the vector field 
$X \in \Gamma (T(T^*Q))$ defined by $i_X \omega_0 = \pi^*\beta$ 
is symplectic with the following well-known properties (see e.~g. 
\cite [Sec 3.2] {BW95}:
\begin {LEMMA} \label {LbetaProjectLem}
Let $\beta \in \Gamma (T^*Q)$ be a closed one-form then the symplectic 
flow $\phi_s$ of the symplectic vector field $X$ defined by 
$i_X \omega_0 = \pi^*\beta$ is complete and given by
\BEQ {XFlow}
    \phi_s (\alpha_q) = \alpha_q - s\beta (q)
\EEQ
for $\alpha_q \in T^*_q Q$ and $s \in \mathbb R$ and for
$L_{s\beta} := \graph (s\beta)$ we have
\BEQ {FlowQLbeta}
    \phi_{-s} (i(Q)) = i_{s\beta} (L_{s\beta}) .
\EEQ
\end {LEMMA}
Clearly $\phi_{-s}$ determines a diffeomorphism 
$\Phi_s: Q \to L_{s\beta}$ by $q \mapsto s\beta(q)$ such that
\BEQ {Flowi}
    \phi_{-s} \circ i = i_{s\beta} \circ \Phi_s
\EEQ
for all $s \in \mathbb R$. For $s = 1$ we obtain from 
(\ref {HamiltonJacobiEq}) and (\ref {FlowQLbeta})
\BEQ {iFlowHE}
    i^* \phi_{-1}^* H = E .
\EEQ

Until now the setting was completely classical and we shall now construct 
a quantum mechanical automorphism as analogue to $\phi_s$ using 
Appendix \ref {TimeApp}: Since $d\pi^*\beta = 0$ the quantum mechanical 
time evolution with respect to this one-form determines a 
one-parameter group of real automorphisms 
$A_s: C^\infty (T^*Q)[[\lambda]] \to C^\infty (T^*Q)[[\lambda]]$
of the star product $\starw$ of the form
\BEQ {AutoAs}
    A_s = \phi^*_s \circ T_s 
    \qquad
    \mbox { where }
    \qquad
    T_s = \id + \sum_{r=1}^\infty \lambda^r T^{(r)}_s
\EEQ
and each $T^{(r)}_s$ is a differential operator which can be computed in 
principle by iterated integrals as in appendix \ref {TimeApp}. Note that 
$A_s \pi^* = \pi^*$. Then the image of the 
two-sided ideal $C^\infty_Q (T^*Q)[[\lambda]]$ under $A_s$ is again a 
two-sided ideal stable under complex conjugation since $A_s$ is real and 
it can easily be determined: let $C^\infty_{L_{s\beta}} (T^*Q)$ be the set 
of those functions $f$ such that $\supp f \cap i_{s\beta} (L_{s\beta})$ 
is compact, then we have the following lemma:
\begin {LEMMA} \label {OmegasLem}
With the notations from above we have for all $s \in \mathbb R$
\BEQ {AsCQ}
    A_s \left(C^\infty_Q (T^*Q)[[\lambda]]\right) = 
    C^\infty_{L_{s\beta}} (T^*Q)[[\lambda]]  
\EEQ    
and hence the $\mathbb C[[\lambda]]$-linear functional
\BEQ {OmegasDef}
    \omega_s := \omega_\mu \circ A_{-s} : \; 
    C^\infty_{L_{s\beta}} (T^*Q) [[\lambda]] \; \to \; \mathbb C[[\lambda]]
\EEQ
is well-defined and positive and for 
$f \in C^\infty_{L_{s\beta}} (T^*Q)[[\lambda]]$ we have
\BEQ {OmegasExpl}
    \omega_s (f) = \int_{L_{s\beta}} i_{s\beta}^*\left( T_{-s} (f)\right)
                   \; \mu_s
                   \qquad 
                   \mbox { where }
                   \qquad
                   \mu_s := {\Phi_s^{-1}}^*\mu .
\EEQ
\end {LEMMA}                   
\begin {PROOF}
Equation (\ref {AsCQ}) follows easily from (\ref {AutoAs}) and 
(\ref {FlowQLbeta}) since $T^{(r)}_s$ is a differential operator. 
Then the well-definedness of $\omega_s$ is obvious 
and the positivity follows from the fact that $A_{-s}$ is a real 
$\starw$-automorphism. Equation (\ref {OmegasExpl}) is a straight forward 
computation using (\ref {Flowi}) and (\ref {AutoAs}).
\end {PROOF}

Now we can apply proposition \ref {AutoGNSProp} and consider the GNS 
representation $\pi_s$ induced by $\omega_s$ on the Hilbert space 
$\mathfrak H_s := C^\infty_{L_{s\beta}} (T^*Q) [[\lambda]] /\mathcal J_s$ 
where $\mathcal J_s$ is the Gel'fand ideal of $\omega_s$.
\begin {LEMMA}
The Gel'fand ideal $\mathcal J_s$ of $\omega_s$ is given by 
$\mathcal J_s = A_s (\mathcal J_\mu)$ and $\mathfrak H_s$ is canonically 
isometric to $C^\infty_0 (L_{s\beta}) [[\lambda]]$ endowed with the 
Hermitian product
\BEQ {HermitiansDef}
    {\SP {\chi} {\varphi}}_s 
    := \int_{L_{s\beta}} \cc \chi \varphi \; \mu_s
    \qquad
    \mbox { where }
    \qquad 
    \chi, \varphi \in C^\infty_0 (L_{s\beta}) [[\lambda]]
\EEQ
via the $\mathbb C[[\lambda]]$-linear unitary map 
\BEQ {HsIso}
    \psi^{(s)}_f \; \mapsto \;  {\Phi_s^{-1}}^* i^* N A_{-s} (f) 
    \quad \textrm{ and its inverse }\quad \chi 
    \; \mapsto \;
    \phi^{(s)}_{\pi^* \Phi_s^* \chi},      
\EEQ
where $f \in C^\infty_{L_{s\beta}}(T^*Q)[[\lambda]]$ and $\chi \in
C^\infty_0 (L_{s\beta})[[\lambda]]$.
Moreover $A_s$ induces due to proposition \ref {AutoGNSProp} a unitary 
map 
$U_s: C^\infty_0 (Q)[[\lambda]] \to C^\infty_0 (L_{s\beta}) [[\lambda]]$
which is given by
\BEQ {UsDef}    
    U_s : \; \chi \mapsto {\Phi_s^{-1}}^* \chi 
    \qquad
    \mbox { for }
    \chi \in C^\infty_0 (Q)[[\lambda]]
\EEQ
and the induced GNS representation $\pi_s$ on 
$C^\infty_0 (L_{s\beta})[[\lambda]]$ is given by
\BEQ {pisDef}
    \pi_s (f) \chi = U_s \wrep (A_{-s} f) U_s^{-1} \chi 
                   = {\Phi_s^{-1}}^* \wrep (A_{-s} f) \Phi_s^* \chi
\EEQ
for all $f \in C^\infty (T^*Q)[[\lambda]]$ and 
$\chi \in C^\infty_0 (L_{s\beta})[[\lambda]]$.                   
\end {LEMMA}    
\begin {PROOF}
It is a straight forward computation to prove that the first map in 
(\ref {HsIso}) is well-defined, bijective and has the inverse as in 
(\ref {HsIso}) due to $A_s \pi^* = \pi^* = N\pi^*$. The fact that 
(\ref {HsIso}) is isometric is computed the same way. Using the 
definition of the unitary map as in proposition \ref {AutoGNSProp} 
and (\ref {HsIso}) the equations (\ref {UsDef}) and (\ref {pisDef}) 
easily follow.
\end {PROOF}

The WKB expansion is now obtained form the following eigenvalue 
problem:
We consider the case $s=1$ then $L_\beta \subseteq H^{-1} (\{E\})$ 
and (\ref {iFlowHE}) is valid. Then we ask for an eigenfunction 
$\tilde \chi \in C^\infty_0 (L_\beta) [[\lambda]]$ of $\pi_1 (H)$ with 
eigenvalue $E$, i.~e.
\BEQ {HEigen}
    \pi_1 (H) \tilde \chi = E \tilde \chi
\EEQ
and due to (\ref {pisDef}) it is equivalent to the corresponding 
eigenvalue problem in $C^\infty_0 (Q)[[\lambda]]$ namely
\BEQ {AHEigen}
    \wrep (A_{-1} H) \chi = E \chi 
\EEQ
where $\chi = \Phi_1^* \tilde \chi$. Now (\ref {AHEigen}) is in each order 
of $\lambda$ a coupled linear partial differential equation for 
$\chi_r \in C^\infty_0 (Q)$ if we write 
$\chi = \sum_{r=0}^\infty \lambda^r \chi_r$. Hence it is well-defined to 
ask for a solution in a distributional sense, i.~e. we ask for  
$\chi \in C^\infty_0 (Q)'[[\lambda]]$ to solve (\ref {AHEigen}) where 
$C^\infty_0 (Q)'$ denotes the space of distributions on $Q$ which is 
obtained as usual as the topological dual of $C^\infty_0 (Q)$ with 
respect to its locally convex topology. The main observation is now 
that we get a linear (first order) partial differential equation for the 
$\chi_r$ which can be solved recursively. The result is a straight 
forward computation expanding (\ref {AHEigen}) in powers of 
$\lambda$ completely analogously to \cite [Theorem 3] {BW97b} 
(where we used a slightly different notation):
\begin {THEOREM} [Formal WKB Expansion]
Let $H \in C^\infty (T^*Q)$ and $L_\beta \subseteq H^{-1} (\{E\})$ be a 
projectable Lagrangean submanifold such that $L_\beta = \graph (\beta)$ 
with $\beta \in \Gamma (T^*Q)$ and let $A_s$ be the one-parameter group of 
$\starw$-automorphisms induced by $\beta$. Then the WKB eigenvalue problem 
(\ref {AHEigen}) for 
$\chi = \sum_{r=0}^\infty \lambda^r \chi_r \in C^\infty_0 (Q)'[[\lambda]]$
is equivalent to the following recursive system of linear first order 
partial differential equations for $\chi$: \\
The homogeneous WKB transport equation for $\chi_0$
\BEQ {WKBHom}
    i^*\{ \im \phi^*_{-1} H, \pi^*\chi_0 \} 
    + i^*\left( \frac{\im\Delta}{2} \phi^*_{-1} H 
                + \phi^*_{-1}T_{-1}^{(1)} H \right) 
    \chi_0 = 0
\EEQ
and the inhomogeneous WKB transport equation for $\chi_r$, $r \ge 1$
\begin {eqnarray}
    \lefteqn {
    i^*\left\{ \im\phi^*_{-1} H, \pi^*\chi_r \right\} 
    + i^*\left( \frac{\im\Delta}{2} \phi^*_{-1} H 
                + \phi^*_{-1}T_{-1}^{(1)} H \right) 
    \chi_r}  \nonumber \\
    & =  &
    - \sum_{{a+b+c+d=r+1 \atop a,b,c,d \ge 0, d < r}}
    i^* M_a \left(\frac{1}{b!} \left(\frac{\im\Delta}{2}\right)^b 
    \phi^*_{-1} T^{(c)}_{-1} H, \pi^*\chi_d \right) \label {WKBInhom}
\end {eqnarray}
where $\Delta$ as in (\ref {DefDelta}) and $M_a$ is the following 
bidifferential operator as in the standard ordered representation, i.~e.
\BEQ {MaDef}
    i^*M_a (f, \pi^*\chi) 
    = \frac{1}{a! \im^a} i^*\left(
    \frac{\partial^a f}{\partial p_{i_1} \cdots \partial_{p_{i_a}}} 
    \right)
    i_s (\partial_{q^{i_1}}) \cdots i_s (\partial_{q^{i_a}}) 
    \frac{1}{a!} D^a \chi 
\EEQ    
and $T^{(c)}_{-1}$ as in (\ref {AutoAs}).
\end {THEOREM}
\begin {PROOF}
The proof is a simple computation of $\wrep (A_{-1} H) \chi$ by expanding 
this in powers of $\lambda$ and (\ref{iFlowHE}).
\end {PROOF}

Note that the question whether this system has a formal solution $\chi$ 
and whether this solution can be `sumed up' to an eventually regular 
distribution after the substitution $\lambda \to \hbar$ is not answered by 
this theorem. Nevertheless it provides a rather explicit recursion scheme 
for the formal eigendistribution $\chi$ 
(the only problem are the operators $T^{(c)}_{-1}$ which are less explicit 
but could be computed by iterated integrals as in appendix \ref {TimeApp}). 
In physical applications the Hamiltonian is often of the form 
(\ref {HfreeV}) and in this case the above recursion scheme is simplified 
even more due to the following lemma
which implies that in the WKB transport equations all terms involving 
the operators $T^{(c)}_{-1}$ with $c \ge 1$ applied to such Hamiltonians 
vanish and thus all operators in the recursion are explicitly given.
\begin {LEMMA}
Let $H \in C^\infty_{pp} (T^*Q)$ be at most quadratic in the momenta and 
let $\beta \in \Gamma (T^*Q)$ be a closed one-form and 
$A_s = \phi_s^* \circ T_s$ the corresponding time development operator. 
Then $A_s H = \phi_s^* H$ and thus $T_s H = H$.
\end {LEMMA}
\begin {PROOF}
This is a simple consequence of the homogeneity and the 
Weyl type property of $\starw$.
\end {PROOF}

Note furthermore that the usual WKB phase can be recovered by the same 
argument as in \cite {BW97b} form the $\starw$-automorphism $A_s$ 
interpreted as the conjugation by the {\em star-exponential} 
$e^{\starw\frac{\im \pi^*S}{\hbar}} = e^{\frac{\im \pi^*S}{\hbar}}$ 
(see \cite {BFFLS78} for a definition) where locally $S$ is a solution of 
the Hamilton Jacobi Equation, i.~e. $dS = \beta$ and $\lambda$ is 
substituted by $\hbar$.

\section{The trace for homogeneous star products on $T^*Q$}
\label {TraceSec}

This section was motivated by a conversation with Markus Pflaum who
proved a particular case of the result we present in this section by
using symbol calculus for pseudo differential operators on Riemannian
manifolds (cf.~\cite{Pfl97b}). Compare also with \cite {CFS92} for the
case of a compact Riemannian configuration space where
another calculus for pseudo differential operators was used.

We shall show in this section in a more algebraically way that the
integration over $T^*Q$ is a trace (i.~e. a linear form vanishing on
star-commutators of two functions where one has compact support) for
{\em all homogeneous}
star products and hence all homogeneous star products are strongly closed
in the sense of \cite {CFS92}. In particular the star products
$\stars$, $\starw$, and the Fedosov star product $\starf$ defined in
\cite{BNW97a} are strongly closed.
\begin {LEMMA} \label {HomDiffOpLem}
Let $D: C^\infty (T^*Q) \to C^\infty (T^*Q)$ be a
homogeneous differential operator of degree $-r$ with 
$r \ge 1$, i.~e. $[\Lie_\xi, D] = -rD$, then
\BEQ {HomDiffOpClosed}
    \int_{T^*Q} D(g) \: \Omega = 0
\EEQ
for all $g \in C^\infty_0 (T^*Q)$ where
$\Omega = \frac{1}{n!} \omega_0^{\wedge n}$ is the symplectic volume
form.
\end {LEMMA}
\begin {PROOF}
First we consider a function $g$ such that $\supp g$ is contained in the
domain of a local bundle chart with coordinates
$q^1, \ldots, q^n, p_1 \ldots, p_n$. Then due to the homogeneity 
$D$ is locally given by
\[
    D(g) = \sum_{I,K \mbox{ with } |I| \ge r} D^K_I
           \frac{\partial^{|I| + |K|} g}{\partial p_I \partial q^K}
\]
where we used the usual notation for multiindices
$I = (i_1, \ldots, i_{|I|})$. The homogeneity of $D$ implies
$\Lie_\xi D^K_I = (|I| - r) D^K_I$ for the coefficient functions
$D_I^K$. Let $l := |I|$ and $I' := (i_1, \ldots, i_{l-1})$ then
\[
    \partial_{p_{i_l}} \left( D^K_I \frac{\partial^{|I|-1+|K|} g}
    {\partial p_{I'} \partial q^K} \right)
    =
    \left( \partial_{p_{i_l}} D^K_I \right)
    \frac{\partial^{|I|-1+|K|} g}
    {\partial p_{I'} \partial q^K}
    +
    D^K_I \frac{\partial^{|I|+|K|} g}
    {\partial p_{I} \partial q^K} 
\]
since $|I| \ge r \ge 1$.
Hence we can conclude by induction on $l \ge r$ since for $l=r$ the
first term vanishes due to the homogeneity of $D^K_I$
that there exist smooth functions $N^K_{iJ}$ defined in the domain of
the chart such that
\[
    D^K_J \frac{\partial^{|I|+|K|} g}
    {\partial p_{I} \partial q^K}
    =
    \partial_{p_i} \left( \sum_{|J| \le l-1} N^K_{iJ}
    \frac{\partial^{|J|+|K|} g}
    {\partial p_{J} \partial q^K}
    \right)
\]
which clearly implies (\ref {HomDiffOpClosed}) for such $g$.
Now let $g \in C^\infty_0 (T^*Q)$ be arbitrary
then using a partition of unity $g$ can be decomposed into
a finite sum of functions each having their support in the domain of
a bundle chart. With the local argument from above the proof is 
completed.
\end {PROOF}
\begin {COROLLARY}
Let $*$ be a homogeneous star product for $T^*Q$ and let
$\widehat T \in C^\infty_{pp,k} (T^*Q)$ and $g \in C^\infty_0 (T^*Q)$
then
\BEQ {intTg}
    \int_{T^*Q} \widehat T * g \: \Omega = \sum_{r=0}^k \lambda^r
    \int_{T^*Q} C_r (\widehat T, g) \: \Omega
\EEQ
and analogously for $\int_{T^*Q} g * \widehat T \: \Omega$ where 
$C_r$ denotes the bidifferential operator of $*$ in order $\lambda^r$.
\end {COROLLARY}
\begin {LEMMA}
Let $*$ be a homogeneous star product for $T^*Q$ and 
$\widehat T \in C^\infty_{pp} (T^*Q)[[\lambda]]$ and 
$g \in C^\infty_0 (T^*Q)[[\lambda]]$ then 
\BEQ {intTraceHom}
    \int_{T^*Q} \left( \widehat T * g - g * \widehat T \right) \: \Omega = 0 .
\EEQ
\end {LEMMA}
\begin {PROOF}
Firstly we can assume that $\supp g$ is contained in the domain of a bundle 
chart and extend the statement afterwards by a partition of unity argument 
as in the proof of lemma \ref {HomDiffOpLem}. Since in such a chart 
any $\widehat T \in C^\infty_{pp} (T^*Q) [\lambda]$ can be written 
locally as a $*$-polynomial in functions at most linear in the momenta 
with coefficient in $\mathbb C[\lambda]$ 
(due to \cite [Prop.~3.7.{\it ii}\/] {BNW97a}) 
it is sufficient to consider such polynomials. Clearly 
$\int_{T^*Q} (\pi^* \chi * g - g * \pi^* \chi) \: \Omega = 0$
for $\chi \in C^\infty (Q)$ and for 
$\widehat X \in C^\infty_{pp, 1} (T^*Q)$ we obtain from (\ref {intTg})
\[
    \int_{T^*Q} (\widehat X * g - g * \widehat X) \: \Omega
    =
    \lambda \int_{T^*Q} (C_1 (\widehat X,g) - C_1 (g, \widehat X)) 
    \: \Omega
    =
    \im\lambda \int_{T^*Q} \{\widehat X, g\} \: \Omega = 0
\]
since the integral over a Poisson bracket vanishes.
Now a simple induction and the $\mathbb C[[\lambda]]$-linearity of 
$\int\ldots$ completes the proof.
\end {PROOF}

Finally we have to extend the statement of this lemma to arbitrary smooth 
functions:
\begin {LEMMA}
Let $D$ be a differential operator with compact support such 
that for any $\widehat T \in C^\infty_{pp} (T^*Q)$ the integral of 
$D(\widehat T) \: \Omega$ over $T^*Q$ vanishes. Then for any function 
$f \in C^\infty (T^*Q)$
\BEQ {intfNull}
    \int_{T^*Q} D(f) \: \Omega = 0 .
\EEQ
\end {LEMMA}
\begin {PROOF}
We shall prove this by induction on the order of $D$. If $D$ is a 
differential operator of order $0$ then it is just a left multiplication 
by a function $D_0$ having compact support. Now if 
$\int_{T^*Q} D_0 \widehat T \: \Omega = 0$ then $D_0 = 0$ since the 
functions polynomial in the momenta are uniformly dense on any 
compactum in all smooth functions due to the Stone-Weierstra{\ss} theorem. 
Now let $D$ be a differential operator of order $k$ and we may again 
firstly assume that $\supp D$ is contained in a local bundle chart and 
extend the result afterwards by a partition of unity. Writing 
$D (f) = \sum_{|I| \le k} D^I \frac{\partial^{|I|} f}{\partial x^I}$ we 
obtain by integration by parts that 
\[
    \int_{T^*Q} D(f) \: \Omega = 
    \int_{T^*Q} \tilde D (f) \: \Omega
    \quad
    \mbox { where }
    \quad
    \tilde D (f) = - \sum_{|I| = k} 
                     \frac{\partial D^I}{\partial x_{i_k}} 
                     \frac{\partial^{|I|-1} f}{\partial x^{I'}}
                     + 
                     \sum_{|I| < k}
                     D^I \frac{\partial^{|I|} f}{\partial x^I} 
\]
is a differential operator of order $k-1$ and obviously 
$\int_{T^*Q} \tilde D (\widehat T) \: \Omega = 0$ for all
$\widehat T \in C^\infty_{pp} (T^*Q)$. Hence by induction 
the proof is complete.
\end {PROOF}

Collecting the results we obtain that any homogeneous star product on 
$T^*Q$ is strongly closed:
\begin {THEOREM} 
Let $*$ be a homogeneous star product for $T^*Q$ and 
$f, g \in C^\infty(T^*Q)[[\lambda]]$ where the coefficients of $g$ have
compact support. Then we have
\BEQ {TraceProperty}
    \int_{T^*Q} (f * g - g * f) \: \Omega = 0
\EEQ
and hence any homogeneous star product is strongly closed,
in particular $\stars$, $\starw$, and $\starf$.
\end {THEOREM}

%
%

\appendix

\section {Formal completed Newton-Puiseux series}
\label {CNPApp}

In this appendix we shall briefly remember the definition and some basic 
facts about formal Laurent series, formal Newton-Puiseux series 
(NP series) and formal completed Newton-Puiseux series (CNP series). For 
proofs and further references we mention \cite {BW96b,Rui93}. These 
formal series will generalize the formal power series in a natural way 
allowing more general exponents of the formal parameter. We define the 
allowed exponents in the following way: a subset $S \subset \mathbb Q$ is 
called L-admissable iff $S$ has a smallest element and 
$S \subset \mathbb Z$, it is called NP-admissable iff it has a smallest 
element and there exists a natural number $N$ such that 
$N\cdot S \subset \mathbb Z$, it is called CNP-admissable iff $S$ has a 
smallest element and $S \cap [i,j]$ is finite for any $i,j \in \mathbb Q$.
Now let $\field K$ be a field, $V$ a vector space over $\field K$ and 
$f : \mathbb Q \to V$ a map, then we define the $\lambda$-support 
$\lsupp f$ of $f$ by
$\lsupp f := \{ q \in \mathbb Q \; | \; f(q) \ne 0 \}$. Then the formal 
Laurent series, the formal NP series and the formal CNP series all with 
coefficients in $V$ are defined by
\BEQ {LNPCNPDef}
    \begin {array} {c}
    \LS V :=  \{ f: \mathbb Q \to V 
              \; | \; \lsupp f \mbox { is L-admissable } \} \\
    \NP V :=  \{ f: \mathbb Q \to V 
              \; | \; \lsupp f \mbox { is NP-admissable } \} \\
    \CNP V := \{ f : \mathbb Q \to V
              \; | \; \lsupp f \mbox { is CNP-admissable } \} .
    \end {array}
\EEQ                                      
Note that clearly $\LS V \subseteq \NP V \subseteq \CNP V$ are all 
$\field K$-vector spaces. Moreover we denote by $V(\lambda)$ the subspace 
of $\LS V$ of those elements with finite support and similarly 
$V \langle \lambda \rangle$. An element $f \in \CNP V$ is now usually 
written as formal series in the formal parameter $\lambda$
\BEQ {fFormalSeries}
    f = \sum_{q \in \lsupp f} \lambda^q f_q
    \quad 
    \mbox { where }
    f_q := f(q) .
\EEQ
Then the vector space structure of $\CNP V$ means just addition of the 
coefficients of the same power of $\lambda$ and scalar multiplication of 
each coefficient and the vector space $V$ itself can be identified 
with the subspace of those elements in $\CNP V$ with exponent $0$ via the 
linear and injective map $V \ni v \mapsto \lambda^0 v \in \CNP V$.
Note that $\CNP V$ and hence all its subspaces can be 
metrisized by the following construction: we define for 
$0 \ne f \in \CNP V$ the order $o(f) := \min (\lsupp f)$ and set 
$o(0) := +\infty$ and define $\varphi (f) := 2^{-o(f)}$ resp. 
$\varphi (0) = 0$ and for $f, g \in \CNP V$ we set
\BEQ {CNPMetricDef}
    d_\varphi (f, g) := \varphi (f - g)
\EEQ
which turns out to be an ultra-metric. Then $V[[\lambda]]$ and $\LS V$ 
are known to be complete metric spaces with respect to $d_\varphi$ and 
$V[\lambda]$ is dense in $V[[\lambda]]$ as well as $V(\lambda)$ is dense 
in $\LS V$. Moreover it was shown in \cite {BW96b} that $\CNP V$ is a 
complete metric space too and $V\langle \lambda \rangle$ as well as 
$\NP V$ are dense subspaces. The topology induced by $d_\varphi$ is 
usually called the $\lambda$-adic topology and the induced topology for 
$V \subseteq \CNP V$ is the discrete topology.

Next we consider an algebra $\mathcal A$ over $\field K$ and define for 
$a, b \in \CNP{\mathcal A}$ a product by
\BEQ {CNPProductDef}
    ab = \left( \sum_{q \in \lsupp a} \lambda^q a_q \right) 
         \left( \sum_{p \in \lsupp b} \lambda^p b_p \right)
    := \sum_{t \in \lsupp a + \lsupp b} \lambda^t
       \sum_{q+p=t} a_q b_p
\EEQ
where 
$\lsupp (ab) = \lsupp a + \lsupp b := \{q+p \; | \; q \in \lsupp a, p \in 
\lsupp b \}$
which turns out to be again a CNP-admissable subset of $\mathbb Q$. 
Note that in each order of $\lambda$ the sum is finite and hence $ab$ is 
again a well-defined element in $\CNP{\mathcal A}$. The following 
proposition is proved as in the case of formal power series:
\begin {PROPOSITION} \label {CNPAlgModulProp}
Let $\mathcal A$ be an algebra over $\field K$ and let $\CNP{\mathcal A}$ 
be endowed with the product (\ref {CNPProductDef}) then $\CNP{\mathcal A}$ 
is again a $\field K$-algebra which is associative resp. commutative 
resp. unital iff $\mathcal A$ is associative resp. commutative resp. 
unital. Moreover 
$\mathcal A[[\lambda]] \subseteq \LS {\mathcal A} 
\subseteq \NP {\mathcal A}$ are subalgebras of $\CNP{\mathcal A}$. Let in 
addition $V$ be an $\mathcal A$-module then $\LS V$ resp. $\NP V$ 
resp. $\CNP V$ is an $\LS{\mathcal A}$- resp. $\NP {\mathcal A}$- resp. 
$\CNP {\mathcal A}$-module with the multiplication analogously to 
(\ref {CNPProductDef}).
\end {PROPOSITION}
In particular one can show that (\ref {CNPProductDef}) applied for the 
field $\field K$ itself defines on $\LS {\field K}$ resp. 
$\NP{\field K}$ resp. $\CNP{\field K}$ again the structure of a field 
and in this case $\varphi$ defines a non-archimedean and non-trivial 
absolute value for these fields. In the case when $\field K$ is 
algebraically closed the Newton-Puiseux theorem ensures that 
$\NP{\field K}$ is again algebraically closed and by K\"{u}rsch\'{a}k's theorem 
one obtains that $\CNP{\field K}$ is algebraically closed too and metric 
complete with respect to $d_\varphi$. These two features were the main 
motivation to consider formal CNP series instead of formal power series or 
formal Laurent series. As application of proposition 
\ref {CNPAlgModulProp} one obtains that $\LS V$ resp. $\NP V$ resp. 
$\CNP V$ are vector spaces over the fields $\LS {\field K}$ resp. 
$\NP{\field K}$ resp. $\CNP{\field K}$.

In a next step we consider linear mappings between vector spaces of formal 
series. Let $V$ and $W$ be vector spaces over $\field K$ and 
$\phi : V \to W$ a $\field K$-linear map. Then there exists a unique 
$\field K [[\lambda]]$-linear map 
$\tilde \phi : V[[\lambda]] \to W[[\lambda]]$ such that 
$\tilde \phi | V = \phi$ which is simply obtained by the 
`$\field K[[\lambda]]$-linear continuation' of $\phi$. The same 
result is true for formal Laurent, NP, and CNP series and for simplicity 
we shall always identify $\phi$ with its continuation and use the 
same symbol. Using this identification we get the following 
natural inclusions:
\BEQ {HomInclusion}
    \begin {array} {c}
        \Hom_{\field K} (V, W) [[\lambda]]
        \subseteq \Hom_{\field K[[\lambda]]} (V[[\lambda]], W[[\lambda]])
        \\
        \LS{\Hom_{\field K} (V, W)}   
        \subseteq \Hom_{\LS{\field K}} (\LS V, \LS W)
        \\
        \NP {\Hom_{\field K} (V, W)} 
        \subseteq \Hom_{\NP{\field K}} (\NP V, \NP W)
        \\
        \CNP {\Hom_{\field K} (V, W)}
        \subseteq \Hom_{\CNP{\field K}} (\CNP V, \CNP W)
    \end {array}
\EEQ        
Similar inclusions are valid for multilinear maps as well.
In the case of formal power series it is known that any 
$\field K[[\lambda]]$-linear map $\phi : V[[\lambda]] \to W[[\lambda]]$ 
is of the form $\phi = \sum_{r=0}^\infty \lambda^r \phi_r$ with 
$\phi_r \in \Hom_{\field K} (V, W)$ and hence the first inclusion in 
(\ref {HomInclusion}) is in fact an equality \cite [Prop. 2.1] {DL88}. 
This is in general no longer the case for the other three inclusions: 
let $\dim V = \infty$ and $\dim W > 0$ and let 
$\{e_n\}_{n \in \mathbb N}$ be 
a set of linear independent vectors in $V$ and hence they are still 
linear independent over $\LS{\field K}$ in $\LS V$. Complete this set 
to a base of $\LS V$ by some vectors $\{f_\nu \}_{\nu \in I}$ and define a 
$\LS {\field K}$-linear map $\phi : \LS V \to \LS W$ 
by $\phi (e_n) := \lambda^{-n} w$ and $\phi (f_\nu) := 0$ where 
$0 \ne w \in W$ is some chosen vector. Then clearly 
$\phi \not\in \LS {\Hom_{\field K} (V, W)}$. 
Nevertheless we are mainly interested in those linear maps which can be 
written as formal series in $\field K$-linear maps.

Now we consider certain continuous functions between vector spaces of 
formal series: let $V, W$ be $\field K$-vector spaces and let 
$T: \CNP V \to \CNP W$ be a not necessarily linear map such that there 
exists a $q \in \mathbb Q$ such that
\BEQ {TraisesDegree}
    o (T(v) - T(v')) \ge o (v-v') + q
\EEQ 
for all $v, v' \in \CNP V$. Then we say that $T$ raises the degree 
at least by $q$. For such maps we have the following property which is 
proved straight forward:
\begin {LEMMA} \label {LipschitzContLem}
Let $T: \CNP V \to \CNP W$ be a map raising the degree at least by 
$q \in \mathbb Q$. Then $T$ is Lipschitz-continuous with respect to the 
metric $d_\varphi$ with Lipschitz-constant $2^{-q}$. In particular 
those $\CNP{\field K}$-linear maps of the form 
$\phi = \sum_{q\in \lsupp \phi} \lambda^q \phi_q 
\in \CNP{\Hom_{\field K} (V, W)}$ are Lipschitz-continuous with 
Lipschitz-constant $2^{-\min{(\lsupp \phi)}}$. 
The same is true for formal power, Laurent, and NP series.
\end {LEMMA}
Then a nice consequence of the metric completeness of $V[[\lambda]]$, 
$\LS V$, and $\CNP V$ is the following formal version of Banach's fixed 
point theorem:
\begin {PROPOSITION} [Formal Banach's fixed point theorem] 
\label {BanachProp}
Let $\field K$ be a field and $V$ a vector space over $\field K$ and
$T: \CNP V \to \CNP V$ a map raising the degree at least by $q > 0$. Then 
there exists a unique fixed point $v_0 \in V$ of $T$, i.~e.
\BEQ {TfixedPoint}
    T(v_0) = v_0
\EEQ
and $v_0$ can be obtained by $v_0 = \lim_{n\to \infty} T^n (v)$ where 
$v \in \CNP V$ is arbitrary and the limit is with respect to the 
$\lambda$-adic topology. The same result holds for formal power series and 
Laurent series.
\end {PROPOSITION}
Note that the proposition is not true in general for the formal NP series 
since $\NP V$ is not Cauchy-complete. The above proposition is useful 
in many situations in deformation quantization where one has to 
`solve an equation by recursion', e.~g. in the construction of 
Fedosov star products or for the construction of the time development 
operator (see e.~g. appendix \ref {TimeApp}).

In the following we shall consider the possibility to extend several 
structures which are defined for formal power series to formal Laurent, NP, 
and CNP series. In particular we are interested in algebra deformations in 
the sense of Gerstenhaber (see e.~g. \cite {GS88}): Let 
$(\mathcal A, \mu_0)$ be an algebra over $\field K$ and let 
$\mu_i : \mathcal A \times \mathcal A \to \mathcal A$ be bilinear maps for 
$i \ge 1$. Then we consider the formal deformation 
$(\mathcal A[[\lambda]], \mu)$ where 
$\mu := \sum_{r=0}^\infty \lambda^r \mu_r$ as a module over 
$\field K[[\lambda]]$. (Usually one has additional conditions for the 
deformation $\mu$, for example it should be associative if $\mu_0$ is 
associtive etc.)
\begin {LEMMA}
\label {ExtendLem}
Let $(\mathcal A, \mu_0)$ be an algebra over $\field K$
and let $(\mathcal A[[\lambda]], \mu = \sum_{r=0}^\infty \lambda^r \mu_r)$ 
be a deformation of $(\mathcal A, \mu_0)$. Then $(\LS{\mathcal A}, \mu)$ 
resp. $(\NP{\mathcal A}, \mu)$ resp. $(\CNP{\mathcal A}, \mu)$ are 
algebras over $\LS{\field K}$ resp. $\NP{\field K}$ resp. 
$\CNP{\field K}$ which are associative resp. Lie resp. commutative iff 
$(\mathcal A[[\lambda]], \mu)$ is associative resp. Lie resp. commutative. 
Let $(\mathcal B[[\lambda]], \tilde \mu)$ be another deformed algebra and 
$\phi: \mathcal A[[\lambda]] \to \mathcal B[[\lambda]]$ a 
$\field K[[\lambda]]$-linear map. Then $\phi$ is an 
(anti-) homomorphism of $\field K[[\lambda]]$-modules iff the extension 
$\phi: (\LS{\mathcal A}, \mu) \to (\LS {\mathcal B}, \tilde \mu)$ resp.
$\phi: (\NP{\mathcal A}, \mu) \to (\NP {\mathcal B}, \tilde \mu)$ resp.
$\phi: (\CNP{\mathcal A}, \mu) \to (\CNP {\mathcal B}, \tilde \mu)$
is an (anti-) homomorphism of $\LS{\field K}$- resp. $\NP{\field K}$-
resp. $\CNP{\field K}$-algebras. The same is true for 
a $\field K[[\lambda]]$-linear derivation $D$ of
$(\mathcal A[[\lambda]], \mu)$.
\end {LEMMA}
This easy lemma is very useful since it justifies to perform many 
calculation only in the setting of formal power series and extend the 
results afterwards.

A last important concept in particular for the general GNS construction as 
proposed in \cite {BW96b} is the notion of positivity. Here we consider an 
ordered field $\field R$ and a quadratic field extension 
$\field C := \field R (\im)$ where $\im^2 := -1$ (Of course we have in 
mind to use the fields of real and complex numbers). 
Then $\field R[[\lambda]]$ 
is an ordered ring and $\LS{\field R}$, $\NP{\field R}$, and
$\CNP{\field R}$ are ordered fields by the following definition: 
an element $a = \sum_{q\in \lsupp a} \lambda^q a_q$ 
is called positive iff $a_{q_0} > 0$ in $\field R$ where 
$q_0 := \min (\lsupp a)$. Complex conjugation in $\field C$ is defined as 
usual and extended to $\field C [[\lambda]]$, $\LS{\field C}$, 
$\NP{\field C}$, and $\CNP{\field C}$ by the definition that the formal 
parameter (and all its powers) should be real: $\cc \lambda := \lambda$. 
Note that the topology induced by the order coincides with the metric 
topology induced by $d_\varphi$.

Now let $(\mathcal A, \mu_0)$ be an associative $\field C$-algebra and 
$\mu = \sum_{r=0}^\infty \lambda^r \mu_r$ an associative deformation of 
$\mu_0$. Let moreover 
$^*: \mathcal A[[\lambda]] \to \mathcal A[[\lambda]]$ be an involutive
$\field C[[\lambda]]$-antilinear anti-automorphism of the form 
$^* = \sum_{r=0}^\infty \lambda^r {}^{*_r} 
\in \Hom_{\field K} (\mathcal A, \mathcal A)[[\lambda]]$. 
Then a $\field C[[\lambda]]$-linear functional 
$\omega: \mathcal A[[\lambda]] \to \field C[[\lambda]]$ is called positive 
iff for all $A \in \mathcal A[[\lambda]]$ we have 
$\omega (\mu(A^*, A)) \ge 0$ in 
$\field R[[\lambda]] \subset \field C[[\lambda]]$ and analogously for 
formal Laurent, NP, and CNP series. Then the following lemma shows that 
the positivity of a linear functional is already determined by the 
elements in $\mathcal A$:
\begin {LEMMA} \label {ExtendPosFunctLem}
Let $\mathcal A, \mu$ and $^*$ be as above and let 
$\omega: \CNP {\mathcal A} \to \CNP {\field C}$ be a 
$\CNP{\field C}$-linear functional of the form 
$\omega = \sum_{q \in \lsupp \omega} \lambda^q \omega_q 
\in \CNP{(\mathcal A^*)}$ (where $\mathcal A^*$ denotes the 
algebraic dual of $\mathcal A$) such that for all 
$A \in \mathcal A$ one has $\omega (\mu(A^*, A)) \ge 0$.
Then $\omega$ is a positive functional on $\CNP{\mathcal A}$.
The analogous result holds for formal power, Laurent, and NP series.
\end {LEMMA}
\begin {PROOF}
Let $A = \sum_{q \in \lsupp A} \lambda^q A_q \in \CNP{\mathcal A}$ then we 
have to prove $\omega (A^*A) \ge 0$. Using the Cauchy-Schwartz inequality
for $\omega$ applied for $A_q \in \mathcal A$ which is valid due to the 
assumption we notice that if $\omega(A_q^*A_q) = 0$ then non of the terms 
involving $A_q$ in $\omega (A^*A)$ contributes. Hence we can assume that 
$\omega (A_q^* A_q) > 0$ for all $q \in \lsupp A$. Now let 
$\lsupp A = \{q_0 < q_1 < \ldots\}$ and we can furthermore assume 
$q_0 = 0$. Let $c := \omega (A_0^* A_0) > 0$, 
$a := \omega (A_{q_1}^* A_{q_1}) > 0$ and $b := \omega (A_0^* A_{q_1})$ 
then the positivity of $\omega$ applied for $A_0$ and $A_{q_1}$ implies 
$\omega (A^*_{q_1} A_0) = \cc b$ and $b\cc b \le ac$. 
This implies on the other hand that 
$\omega ((A_0 + tA_{q_1})^* (A_0 + tA_{q_1})) = at^2 + (b+\cc b) t + c$ 
is non-negative for all $t = \cc t \in \CNP{\field R}$ and hence it is 
non-negative for $t = \lambda^{q_1}$. Now one proceeds analogously by 
induction to obtain that $\omega$ is positive for all finite sums 
$A_N := \sum_{q \le N} \lambda^q A_q$ but then the continuity of $\mu$, 
$^*$ and $\omega$ according to lemma \ref {LipschitzContLem} guarantees 
that $\omega$ is in fact a positive linear functional on all 
$A \in \CNP{\mathcal A}$ since the order topology for $\CNP{\field R}$ 
coincides with the metric topology induced by $d_\varphi$.
\end {PROOF}

\section {Time development in deformation quantization}
\label {TimeApp}

In this appendix we shall briefly remember some well-known facts 
about the time development in deformation quantization 
(see e.~g. \cite [Sec.~5] {BW96b} and Fedosov's book 
\cite [Sec.~5.4] {Fed96}). Let 
$(M, \omega)$ be a symplectic manifold and let $*$ be a star 
product for $C^\infty (M)[[\lambda]]$. Let $X$ be a symplectic vector 
field with complete flow $\phi_t$ and $\beta = i_X \omega$ the 
corresponding closed but not necessarily exact one-form.
Then locally $\beta$ is exact and we have $\beta = dH$. We use this 
locally defined Hamiltonian $H$ to define star product commutators 
locally and notice that $\ad (H)$ does not depend on the choice of $H$ 
but only on $\beta$ and hence one obtains a globally defined map  
which we shall denote by $\ad (\beta)$. Then we consider the 
Heisenberg equation of motion
\BEQ {HeisenbergEq}
    \frac{d}{dt} f(t) = \frac{\im}{\lambda} \ad (\beta) f(t)
\EEQ
with respect to $\beta$ and ask for a solution 
$t \mapsto f(t) \in C^\infty (M)[[\lambda]]$ for a given initial 
value $f(0)$:
\begin {THEOREM} \label {TimeDevelUniqueTheo}
Let $X$ be a symplectic vector field with complete flow and let
$\beta = i_X \omega$. Then the Heisenberg equation of motion 
(\ref{HeisenbergEq}) has a unique solution $f(t)$ for any initial 
value $f(0) \in C^\infty (M)[[\lambda]]$ and $t \in \mathbb R$.
\end {THEOREM}
\begin {PROOF}
Equation (\ref {HeisenbergEq}) is equivalent to the equation 
$\frac{d}{dt} g(t) = \phi^*_{-t} \circ \hat H \circ \phi^*_t g(t)$
where $\hat H := \frac{\im}{\lambda} \ad (\beta) - \Lie_X $ and 
$g(t) = \phi^*_{-t} f(t)$. This equation can be rewritten as 
a fixed point equation by integrating over $t$. Then 
proposition \ref {BanachProp} shows that there exists a unique 
solution since $\hat H$ raises the $\lambda$-degree at least by one.
\end {PROOF}

This theorem allows us to define 
a quantum mechanical time development operator $A_t$ 
for the symplectic vector field $X$ such that 
$f(t) = A_t f$ is the unique solution with initial condition $f$. 
Then $A_t$ is a $\mathbb C[[\lambda]]$-linear map for all 
$t \in \mathbb R$ and clearly $A_0 = \id$. Moreover $A_t$ 
clearly satisfies the Heisenberg equation as operator equation 
\BEQ {AtHeisenberg}
    \frac{d}{dt} A_t = \frac{\im}{\lambda} \ad (\beta) A_t .
\EEQ
The quantum mechanical time development operator $A_t$ 
is obtained by a quantum correction of the classical 
time development operator which is just $\phi^*_t$. 
We denote this correction by
\BEQ {TtDef}
    T_t := \phi^*_{-t} \circ A_t 
\EEQ
and prove the following proposition by a straight forward 
computation:
\begin {PROPOSITION}
The operator $T_t$ is a formal power series of differential 
operators $T_t = \id + \sum_{r=1}^\infty \lambda^r T^{(r)}_t$ and 
$T_t$ satisfies the following differential equation with initial 
condition $T_0 = \id$
\BEQ {TtHeisenberg}
    \frac{d}{dt} T_t = \phi^*_{-t} \circ \hat H \circ \phi^*_t \circ T_t
\EEQ
and the equivalent integral equation
\BEQ {TtIntegral}
    T_t = \id + \int_0^t \phi^*_{-\tau} \circ \hat H \circ 
    \phi^*_\tau \circ T_\tau d\tau
\EEQ
where $\hat H = \frac{\im}{\lambda} \ad (\beta) - \Lie_X$ 
is defined as in the proof of theorem \ref {TimeDevelUniqueTheo}. 
If the star product is of the Vey type then $T^{(r)}_t$ is a 
differential operator of order $2r$.
\end {PROPOSITION}
In a last step we prove that $A_t$ is a one-parameter group of 
automorphisms of the star product using the fact that the solution 
$f(t)$ of (\ref {HeisenbergEq}) is uniquely determined by $f(0)$:
\begin {THEOREM}
The quantum mechanical time development operator $A_t$ of the 
symplectic vector field $X$ with complete classical flow
has the following properties where $\beta = i_X \omega$:
\begin {enumerate}
\item $A_t A_s = A_{t+s} = A_s A_t$ and $A_0 = \id$ for all 
      $t, s \in \mathbb R$.
\item $A_t \ad (\beta) = \ad (\beta) A_t$ for all 
      $t \in \mathbb R$.      
\item $A_t (f * g) = A_t f * A_t g$ for all 
      $f, g \in C^\infty (M)[[\lambda]]$ and $t \in \mathbb R$.
\item $A_{-t}$ is the time development operator for the vector 
      field $-X$ and 
      $(T_t)^{-1} = \phi^*_{-t} \circ T_{-t} \circ \phi^*_t$.
\item If in addition $\cc{f*g} = \cc g * \cc f$ then $A_t$ is a 
      real automorphism, i. e. $\cc{A_t f} = A_t \cc f$.
\end {enumerate}
\end {THEOREM}
Remark: All statements are still correct if one replaces the 
closed one-form $\beta$ by 
$\beta + \sum_{r=1}^\infty \lambda^r \beta_r$ where $\beta_r$ 
are again closed one-forms which encloses the case of some 
`quantum corrections'. Furthermore all statements can be transfered to the 
case of formal Laurent and CNP series including the corresponding quantum 
corrections (but not necessarily to formal NP series since we used the 
formal Banach's fixed point theorem).

\section*{Acknowledgements}

We would like to thank the old Babylonians for their trick to solve 
quadratic equations and Markus Pflaum for a fruitful discussion motivating
us to add Section 8 of this paper.

\end {document}